\documentclass[11pt]{amsart}

\usepackage{amsmath}
\usepackage{amssymb}
\usepackage{graphicx}

\newtheorem{theorem}{Theorem}[section]
\newtheorem{corollary}[theorem]{Corollary}
\newtheorem{lemma}[theorem]{Lemma}
\newtheorem{proposition}[theorem]{Proposition}
\theoremstyle{definition}

\newtheorem{remark}[theorem]{Remark}

\newcommand{\R}{\mathbb{R}}
\newcommand{\N}{\mathbb{N}}
\newcommand{\D}{\displaystyle}
\newcommand{\norm}[1]{\|#1\|}

\numberwithin{equation}{section}

\title[Nonhomogeneous  equations  involving  Hardy  potentials ]{Solutions of  nonhomogeneous  equations  involving
    Hardy  potentials with singularities on the  boundary }

\author[H. Chen]{Huyuan Chen}

\address[H. Chen]{Department of Mathematics, Jiangxi Normal University,
Nanchang, Jiangxi 330022, PR China}
\email{\tt chenhuyuan@yeah.net}

\author[A. Quaas]{Alexander Quaas}
\address[A. Quaas]{Departamento de Matem\'{a}tica, Universidad T\'{e}cnica Federico,
Santa Mar\'{i}a Casilla V-110, Avda. Espa\~{n}a 1680,
Valpara\'{i}so, Chile}
\email{\tt alexander.quaas@usm.cl}

\author[F. Zhou]{Feng Zhou}
\address[F. Zhou]{Center for PDEs, School of Mathematical Sciences, East China Normal University,
Shanghai Key Laboratory of PMMP,
 Shanghai, 200241, PR China}
\email{\tt fzhou@math.ecnu.edu.cn}



\begin{document}

\begin{abstract}
In this paper,  we present a new  distributional identity  for the  solutions of elliptic equations involving
Hardy potentials with singularities located on the boundary of the domain.
Then we use it to obtain the boundary isolated singular solutions of nonhomogeneous problems.
\end{abstract}

\maketitle


\section{Introduction}

The classical Hardy inequality is stated as following: For any smooth bounded domain $\mathcal{O}$ in $\R^N$ containing the origin, there holds
\begin{equation}\label{0.1}
\int_{\mathcal{O}} |\nabla u|^2 dx \ge c_{N}  \int_{\mathcal{O}}   |x|^{-2} |u|^2 dx ,\quad \forall\,u\in H_0^1(\mathcal{O}),
\end{equation}
with the best constant $c_N=\frac{(N-2)^2}{4}$. The qualitative properties of Hardy inequality and its improved versions have been studied extensively, see for example \cite{ACR,BM,Fi,GP}, motivated by great applications in  the study of  stability of solutions to semilinear elliptic and parabolic equations (cf. \cite{BV1,CM1,Da, PV,V}). 
The isolated singular solutions of Hardy problem with absorption nonlinearity have been studied in \cite{CC,C,GV} and the one with source nonlinearity has been done in \cite{BDT,F}. The related semilinear elliptic problem involving the inverse square potential has been studied  by   variational methods in \cite{D,DD,FF} and the references therein.
In a very recent work \cite{CQZ}, we established a new distributional identity with respect to a specific weighted measure and we then classify the classical isolated singular solutions of
$$-\Delta u+\frac{\mu}{|x|^2}u =f\quad{\rm in}\quad \mathcal{O}\setminus\{0\}, $$
subject to the homogeneous Dirichlet boundary condition with $\mu \geq -c_N$. These results allow us to draw a complete picture of the existence, non-existence and the singularities
for classical solutions for the above problems (cf. \cite{CZ}).


It is of interest to consider the corresponding problem involving Hardy potential with  singularity on the boundary. While the sharp constant $c_N$ in Hardy inequality (\ref{0.1}) could be replaced by $\frac{N^2}{4}$ when the origin is addressed on the boundary of the domain, see \cite[Corollary 2.4]{TFT}, also \cite{Ca,Ca1,FM2}.

Let $\Omega$ be a smooth bounded domain in $\R^N$ with $0 \in \partial\Omega$. We study boundary isolated singular solutions of
nonhomogeneous problems:
 \begin{equation}\label{eq 1.1fg}
 \arraycolsep=1pt\left\{
\begin{array}{lll}
 \displaystyle   \mathcal{L}_\beta u= f\quad
   {\rm in}\;\;  \Omega,\\[1.5mm]
 \phantom{   L_\beta   }
 \displaystyle  u=g\quad  {\rm   on}\;\; \partial{\Omega}\setminus \{0\},
 \end{array}\right.
\end{equation}
where $f\in C^\gamma_{loc}(\bar\Omega\setminus\{0\})$ with $\gamma\in(0,1)$, $g\in C(\partial\Omega\setminus\{0\})$ and
$\mathcal{L}_\beta:= -\Delta   +\frac{\beta}{|x|^{2 }}$
is the Hardy operator which is singular at $0$ (with $N\ge 2$, $\beta\ge \beta_0:=-\frac{N^2}{4}$).
    Recall that for $\beta\ge\beta_0$,  the problem
\begin{equation}\label{eq 1.1}
\left\{\arraycolsep=1pt
\begin{array}{lll}
\mathcal{L}_\beta u= 0\quad {\rm in}\;\;  \R^N_+, \\[1mm]
 \phantom{ \mathcal{L}_\beta  }
  u=0\quad {\rm on}\;\;  \partial\R^N_+\setminus \{0\}
  \end{array}
 \right.
\end{equation}
has two special  solutions with the explicit formulas as
\begin{equation}\label{1.1}
 \Lambda_\beta(x)=\left\{\arraycolsep=1pt
\begin{array}{lll}
 x_N|x|^{\tau_-(\beta)}\quad
   &{\rm if}\;\; \beta>\beta_0,\\[1mm]
 \phantom{   }
-x_N|x|^{\tau_-(\beta)}\ln|x| \quad  &{\rm   if}\;\; \beta=\beta_0
 \end{array}
 \right.\quad \;\;{\rm and}\quad \lambda_\beta(x)=x_N|x|^{\tau_+(\beta)},
\end{equation}
where $x=(x',x_N)\in\R^N_+:=\R^{N-1}\times(0,+\infty)$, and
\begin{equation}\label{tau}
 \tau_-(\beta)=-\frac{N}2-\sqrt{ \beta-\beta_0}\quad{\rm and}\quad  \tau_+(\beta)=-\frac{ N}2+\sqrt{ \beta-\beta_0},
\end{equation}
are two roots of $\beta-\tau(\tau+N)=0$.

As in \cite{CZ,CQZ}, we first find a certain distributional identity which shows that the singularity of  solution $\Lambda_\beta$ for (\ref{eq 1.1})
is associated to a Dirac mass.
 Let $C^{1.1}_0( \R^N_+)$ be the set of functions in $C^{1.1}(\overline{\R^N_+})$ vanishing on the boundary and having compact support in $\overline{\R^N_+}$. Then we have
 \begin{theorem}\label{teo 0}
Let  $d\gamma_\beta :=\lambda_\beta(x) dx$ and
\begin{equation}\label{L}
 \mathcal{L}^*_\beta :=-\Delta -\frac{2\tau_+(\beta) }{|x|^2}\,x\cdot\nabla-\frac2{x_N}\frac{\partial}{\partial x_N}, \quad x=(x',x_N) \in\R^N_+.
\end{equation}
Then there holds
\begin{equation}\label{1.2}
 \int_{\R^N_+}\Lambda_\beta   \mathcal{L}^*_\beta(\frac{\zeta}{x_N})\,  d\gamma_\beta  =c_\beta  \frac{\partial \zeta}{\partial x_N}(0),\quad \forall\, \zeta\in C^{1.1}_0( \R^N_+),
\end{equation}
where 
  \begin{equation}\label{cmu}
 c_\beta=\left\{\arraycolsep=1pt
\begin{array}{lll}
 \sqrt{\beta-\beta_0}\,|\mathcal{S}^{N-1}|/N\quad
   &{\rm if}\;\; \beta>\beta_0,\\[1.5mm]
 \phantom{   }
|\mathcal{S}^{N-1}|/N \quad  &{\rm  if}\;\; \beta=\beta_0,
 \end{array}
 \right.
 \end{equation}
and $\mathcal{S}^{N-1}$ is the unit sphere of $\R^N$ and $|\mathcal{S}^{N-1}|$ denotes its $(N-1)$-dimensional Hausdorff measure.
 \end{theorem}

From the distributional identity (\ref{1.2}), $\Lambda_\beta$  is called as   a fundamental solution of (\ref{eq 1.1}).   We remark that when $\beta=0$, $\mathcal{L}^*_0=-\Delta-\frac2{x_N}\frac{\partial}{\partial x_N}$, $\lambda_{\beta}(x)=x_N$ and  (\ref{1.2}) could be reduced to
\begin{eqnarray*}
c_0 \frac{\partial \zeta}{\partial x_N}(0)= \int_{\R^N_+}\Lambda_0   \mathcal{L}^*_0(\frac{\zeta}{x_N})\,  d\gamma_\beta  =  \int_{\R^N_+}\Lambda_0  (-\Delta \zeta) \,  dx,\quad \forall \;\zeta\in C^{1.1}_0(\R^N_+),
\end{eqnarray*}
which coincides with the classical distributional identity proposed in \cite{GV1}.
On this classical subject, it has been vastly expanded in the works \cite{BV,MV1,MV2,MV3,MV4}.

For simplicity, here and in the sequel, we always assume that $\Omega$ is a  bounded $C^2-$ domain satisfying that
\begin{equation}\label{O}
 B_{r_0}^+(0) \subset \Omega\subset B_{R_0}^+(0),
\end{equation}
for some $0<r_0<R_0<+\infty$ where $B_{r}^+(0):= B_r(0)\cap \R^N_+$.
Let $d\omega_\beta(x) := |x|^{\tau_+(\beta)}d\omega(x)$,  where $\omega$ is the Hausdorff measure of $\partial\Omega$.
We can state our main result as follows

\begin{theorem}\label{teo 1}
Let  $\mathcal{L}^*_\beta$  be given by (\ref{L}),  $f\in C^{\theta}_{loc}( \bar\Omega\setminus\{0\} ) $  with $\theta \in(0,1)$, $g\in C(\partial\Omega\setminus\{0\})$.

$(i)$ If  \begin{equation}\label{1.3fg+}
  \int_{\Omega } |f| \,d\gamma_\beta +\int_{\partial\Omega} |g| \,d\omega_\beta<+\infty,
  \end{equation}
then for any $k\in\R$,    problem (\ref{eq 1.1fg}) 
admits a  unique solution $u_k\in C^2(\Omega)\cap L^1(\Omega,|x|^{-1}d\gamma_\beta)$ such that
 \begin{equation}\label{1.2fg}
 \int_{\Omega}u_k  \mathcal{L}_\beta^*(\frac{\xi}{x_N})\, d\gamma_\beta  = \int_{\Omega}   \frac{f\xi}{x_N} \, d\gamma_\beta-\int_{\partial\Omega} g\frac{\partial \xi}{\partial \nu} d\omega_\beta +c_\beta k \frac{\partial \xi}{\partial x_N}(0) ,\quad\forall\, \xi\in   C^{1.1}_0(\Omega),
\end{equation}
where  $\nu$ is the unit outward vector on $\partial \Omega$.

$(ii)$ If $f, g$ are nonnegative and
 \begin{equation}\label{1.3fg-}
 \lim_{r\to0^+} \Big(\int_{\Omega\setminus B_r(0)} f \,d \gamma_\beta +\int_{\partial\Omega\setminus B_r(0)} g \,d\omega_\beta\Big)=+\infty,
  \end{equation}
 then problem (\ref{eq 1.1fg}) has no nonnegative solution.
\end{theorem}

When $g=0$ on $\partial\Omega$ and $f=0$ in $\Omega$, we prove in Proposition \ref{pr 2.2} in Section 3 that problem (\ref{eq 1.1fg}) admits an isolated singular solution $\Lambda^\Omega_\beta$, which has the asymptotic behavior at the origin as
the fundamental function $\Lambda_\beta$. More precisely, we have
\begin{equation}\label{1.3-00}
 \lim_{t\to0^+}\sup_{z\in S^{N-1}_+}\Big(\frac{\Lambda^\Omega_\beta(tz)}{\Lambda_\beta(tz)} -1\Big)=0.
\end{equation}
When $g=0$ on $\partial\Omega$ and $f\in C^{\theta}_{loc}( \bar\Omega\setminus\{0\} ) \cap L^1(\Omega,d\gamma_\beta)$, Theorem \ref{teo 3.1} in Section 4 shows  that problem (\ref{eq 1.1fg})
has a solution $u_f$ verifying the isolated singularity (see Remark \ref{re 3.2})
\begin{equation}\label{1.3}
 \lim_{t\to0^+}\inf_{z\in S^{N-1}_+}\frac{u_{f}(tz)}{\Lambda_\beta(tz)}=0,
\end{equation}
which is less precise than (\ref{1.3-00}) due to the lack of estimates of Green kernel of Hardy operator with singularity on the boundary.
However, when $f=0$ and $g\not=0$, it is not convenient to use (\ref{1.3}) to describe  the singularity of the solution $u_g$, so we may distinguish this by  the distributional identity
$$\int_{\Omega}u_g \mathcal{L}_\beta^*(\frac{\xi}{x_N})\, d\gamma_\beta  = -\int_{\partial\Omega} g\frac{\partial \xi}{\partial \nu} d\omega_\beta,\quad\forall\, \xi\in   C^{1.1}_0(\Omega),$$
All in all, the solution $u_k$ of (\ref{eq 1.1fg}) can be decomposed into three components $k\Lambda^\Omega_\beta$, $u_f$ and $u_{g}$.

The method we use to prove the existence  of solutions for problem (\ref{eq 1.1fg}) is different from  the classical method of the boundary data problem used by
Gmira-V\'{e}ron in \cite{GV1} due to the appearance of Hardy potential. They obtained the very weak solutions by approximating the  Dirac mass at  boundary. Then they considered the limit of the solutions to the corresponding
problem where the convergence is guaranteed by the Poisson kernel. In this paper, we prove the existence of moderate singular solution by using the function $\Lambda_\beta$ to construct suitable solutions of  problem (\ref{eq 1.1fg}) with the zero Dirichlet boundary condition. While for nonzero Dirichlet boundary condition, we transform the boundary data into nonhomogeneous term. However, for $\beta>0$, that transformation can not totally solve  (\ref{eq 1.1fg})  with the  nonzero Dirichlet boundary condition, and our idea is to cut off the boundary data and approximate the solutions.

The rest of the paper is organized as follows. In Section 2, we  start from a comparison principle for $\mathcal{L}_\beta$ and show the moderate singular solution of  (\ref{eq 1.1fg}) when $g=0$.    Section 3 is devoted to prove the distributional identity (\ref{1.2})  for the fundamental solution $\Lambda_\beta$ in $\R^N_+$, to consider its trace, the corresponding   distributional identity in bounded smooth domain.    Section 4 is to study the qualitative properties of the solutions for problem (\ref{eq 1.1fg}) when $g=0$  and  then we give the  proof of Theorem \ref{teo 1} in the case of nonzero boundary data in Section 5.  In what follows,  we denote by $c_i$ a generic  positive constant in the proofs of the results.

\setcounter{equation}{0}
\section{Preliminary}

\subsection{Comparison principle}
We start the analysis from a comparison principle for $\mathcal{L}_\beta$.   Let $\eta_0:[0,+\infty)\to [0,\,1]$ be a decreasing $C^\infty$ function such that
\begin{equation}\label{eta}
 \eta_0=1\quad{\rm in}\;\; [0,1]\quad{\rm and}\quad \eta_0=0\quad {\rm in}\;\; [2,+\infty).
\end{equation}

\begin{lemma}\label{lm ccp}
Let $\Omega$ be a bounded open set in $\R^N_+$, $L: \Omega\times [0,+\infty)\to[0,+\infty)$ be a continuous function satisfying that for any $x\in  \Omega$,
$$L(x,s_1)\ge L(x,s_2)\quad {\rm if}\;\; s_1\ge s_2,$$ then $\mathcal{L}_\beta+L$ with $\beta\ge \beta_0$ verifies the comparison principle,
that is, if $u,\,v\in C^{1,1}(\Omega)\cap C(\bar \Omega)$ verify that
$$\mathcal{L}_\beta u+ L(x,u)\ge \mathcal{L}_\beta v+ L(x,v) \quad {\rm in}\;\;  \Omega
\quad{\rm and}\quad   u\ge  v\quad {\rm on}\;\; \partial\Omega,$$
then
$u\ge v\;{\rm in}\;  \Omega.$

\end{lemma}
 \noindent {\bf Proof.} Let $w=u-v$ and then $w\ge 0$ on  $\partial\Omega$. Denote $w_-=\min\{w,0\}$,  and we claim that $w_-\equiv 0$.
Indeed, if $ \Omega_-:=\{x\in \Omega:\, w(x)<0\}$ is not empty, then it is a bounded
$C^{1,1}$ domain in $\Omega$ and $w_-=0$ on $\partial\Omega$. We observe that
$\Omega_-\subset\R^N_+$ and  then
from Hardy inequality \cite[(1.7)]{Ca} (see also \cite{MM}), it holds that
\begin{eqnarray*}
 0  &=& \int_{ \Omega_-}(-\Delta  w_- +\frac{\beta}{|x|^{2 }}  w_- ) w_- dx+\int_{ \Omega_-}[L(x,u)-L(x,v)]w_-\,dx \\ &\ge&  \int_{ \Omega_-} \left(|\nabla w_-|^2 +\frac{\beta}{|x|^{2 }}  w_-^2\right) dx
 \ge  c_1 \int_{ \Omega_-}w_-^2 dx,
\end{eqnarray*}
then  $w_-=0$ in $\Omega_-$, by the continuity of $w_-$, which is impossible with the definition of $\Omega_-$.\hfill$\Box$\smallskip

 \begin{lemma}\label{lm cp}
Assume that $\beta\ge\beta_0$,     $f_1$, $f_2$ are two functions in $C^{\theta}_{loc} ( \Omega)$ with $\theta\in(0,1)$, $g_1$, $g_2$ are two continuous  functions on $\partial \Omega\setminus\{0\}$, and
$$ f_1\ge f_2\quad {\rm in}\;\; \Omega \quad{\rm and}\quad  g_1\ge g_2\quad {\rm on}\;\; \partial \Omega\setminus\{0\}.$$
Let $u_i$ ($i=1,2$)  be the classical solutions of
$$
 \arraycolsep=1pt\left\{
\begin{array}{lll}
 \displaystyle \mathcal{L}_\beta u = f_i\quad
   {\rm in}\;\;  {\Omega},\\[1.5mm]
 \phantom{ L_\beta     }
 \displaystyle  u= g_i\quad  {\rm   on}\;\; \partial{\Omega}\setminus \{0\}.
 \end{array}\right.
$$
If
\begin{equation}\label{cm}
\lim_{r\to0^+}\inf_{x\in \partial_+ B_r(0)}[u_1(x)-u_2(x)]\Lambda_\beta^{-1}(x) \ge 0,
\end{equation}
where
$\partial_+ B_r(0)=\partial B_r(0)\cap\Omega$.
Then
$u_1\ge u_2\quad {\rm in}\;\; \overline{\Omega}\setminus\{0\}.$

\end{lemma}
 \noindent {\bf Proof.}  Let $w=u_2-u_1$, then $w$ satisfies
 $$\arraycolsep=1pt\left\{
\begin{array}{lll}
 \displaystyle\qquad \mathcal{L}_\beta w \le 0\quad
   {\rm in}\;\;  {\Omega},\\[1.5mm]
 \phantom{ L_\beta     }
 \displaystyle \qquad  w\le 0 \quad  {\rm   on}\;\; \partial{\Omega}\setminus \{0\},\\[1.5mm]
 \phantom{   }
  \displaystyle \lim_{r\to0^+}\sup_{x\in \partial_+ B_r(0)  }w(x)\Lambda_\beta^{-1}(x) \le 0.
 \end{array}\right.$$
 Thus for any $\epsilon>0$, there exists $r_\epsilon>0$ converging to zero as $\epsilon\to0$ such that
 $$w\le \epsilon \Lambda_\beta\quad{\rm on}\quad \partial B_{r_\epsilon}(0)\cap \Omega.$$
We observe that $w\le 0<\epsilon \Lambda_\beta \ {\rm on}\  \partial\Omega\setminus B_{r_\epsilon}(0),$
which implies by Lemma \ref{lm ccp} that
$$w\le \epsilon \Lambda_\beta\quad{\rm in}\quad \overline{\Omega}\setminus\{0\}. $$
Therefore we obtain that $w\le 0$   in $\overline{\Omega}\setminus\{0\}$ which ends the proof.
\hfill$\Box$\medskip

For any $\varepsilon>0$, denote
\begin{equation}\label{lepsilon}
 \mathcal{L}_{\beta,\varepsilon} =-\Delta +\frac{\beta}{|x|^2+\varepsilon}.
\end{equation}
We remark that $\mathcal{L}_{\beta,\varepsilon}$ is strictly elliptic operator and we have the following existence result for related nonhomogeneous problem.
\begin{lemma}\label{lm 2.1}
Assume that $\varepsilon\in(0,\,1)$, $\beta\ge \beta_0$, $\mathcal{L}_{\beta,\varepsilon}$ is given by (\ref{lepsilon}) and
 $f\in C^\theta_{loc}(\Omega)\cap C(\bar \Omega)$ with $\theta\in(0,1)$ and $g\in C(\partial\Omega)$.
Then  the problem
\begin{equation}\label{eq2.2}
 \arraycolsep=1pt\left\{
\begin{array}{lll}
 \displaystyle  \mathcal{L}_{\beta,\varepsilon} u= f \quad
   &{\rm in}\quad  \Omega,\\[2mm]
 \phantom{  \mathcal{L}_{\beta,\varepsilon} }
 u=g\quad  &{\rm   on}\quad \partial  \Omega
 \end{array}\right.
\end{equation}
has a unique classical  solution $u_\varepsilon\in C^2(\Omega)\cap C(\bar \Omega)$, which verifies that
\begin{equation}\label{11.1}
\int_\Omega  u_\varepsilon  \mathcal{L}_\beta^* (\frac{\xi}{x_N})\, d\gamma_\beta=\int_\Omega  \frac{ f\xi}{x_N}   d\gamma_\beta-\int_{\partial\Omega}  g \frac{\partial \xi}{\partial \nu}\, d\omega_\beta
+  \beta\varepsilon \int_\Omega\frac{  u_\varepsilon \xi}{(|x|^2+\varepsilon)|x|^2x_N} \, d\gamma_\beta,
\end{equation}
for any $\xi\in C^{1.1}_0(\Omega)$.

Assume more that $f\ge0$ in $\Omega$ and $g\ge0$ on $\partial\Omega$.
Then the mapping $\varepsilon\mapsto u_\varepsilon$ is decreasing  if $\beta>0$, and is increasing if  $\beta_0\le \beta<0$.

\end{lemma}
 \noindent {\bf Proof.} We first prove the existence of solution to problem (\ref{eq2.2}). We introduce   Poisson kernel $P_\Omega$ of $-\Delta$ in $\Omega$,  and denote  Poisson operator as
$$
\mathbb{P}_\Omega[g](x)=\int_{\partial \Omega} P_\Omega(x,y)g(y)dy.$$
We observe that
$$  \mathcal{L}_{\beta,\varepsilon} \mathbb{P}_\Omega[g]=\frac{\beta}{|x|^2+\varepsilon} \mathbb{P}_\Omega[g] \in C^1(\Omega)\cap C(\bar\Omega).$$
Then the solution of (\ref{eq2.2}) denoted by $u_\varepsilon$,  could be reduced to $u_\varepsilon=\mathbb{P}_\Omega[g]+u_f$, where $u_f$ is the solution of
\begin{equation}\label{eq2.2.1}
 \arraycolsep=1pt\left\{
\begin{array}{lll}
 \displaystyle  \mathcal{L}_{\beta,\varepsilon} u= f-\frac{\beta}{|x|^2+\varepsilon} \mathbb{P}_\Omega[g] \quad
   &{\rm in}\quad  \Omega,\\[2mm]
 \phantom{  \mathcal{L}_{\beta,\varepsilon} }
 u=0\quad  &{\rm   on}\quad \partial  \Omega.
 \end{array}\right.
\end{equation}
For $\beta\ge\beta_0$, a  solution $u_f$ in $H^1_0(\Omega)$ of (\ref{eq2.2.1}) could be derived
by Ekeland's variational methods as the critical point of the functional
$$I(u)=\int_\Omega |\nabla u|^2dx+\beta\int_\Omega \frac{u^2}{|x|^2+\varepsilon} dx-\int_\Omega\Big(f-\frac{\beta}{|x|^2+\varepsilon} \mathbb{P}_\Omega[g]\Big)u dx. $$
That is well-defined in $H^1_0(\Omega)$ since
$ \beta\in(\beta_0,0)$. From the Hardy's inequality in \cite{FM2},  we have that, for any $u\in C_0^2(\Omega)$,
$$\int_\Omega |\nabla u|^2 dx +\beta\int_\Omega \frac{u^2}{|x|^{2}+\varepsilon}dx\ge (\beta-\beta_0)  \int_\Omega |\nabla u|^2 dx,$$
 for $ \beta=\beta_0$, from the improved Hardy inequality in \cite{FM2}, it holds
\begin{eqnarray*}
c_2\int_\Omega  u^2 dx &\le & \int_\Omega |\nabla u|^2 dx-|\beta_0|\int_\Omega \frac{u^2}{|x|^{2}} dx\\
    &<& \int_\Omega |\nabla u|^2 dx-|\beta_0|\int_\Omega \frac{u^2}{|x|^{2}+\varepsilon } dx.
\end{eqnarray*}
Finally it is trivial for the case $ \beta\ge 0$.

 By the standard regularity result (e.g. \cite{HL}), we have that $u_f$ is a classical solution of (\ref{eq2.2.1}).
 Then  problem (\ref{eq2.2}) admits a classical solution  and the uniqueness follows by comparison principle.\smallskip

Finally, we prove (\ref{11.1}).  Multiple $\frac{\lambda_\beta\xi}{x_N}$ with $\xi\in C^{1.1}_0(\Omega)$ and integrate over $\Omega$,  we have that
 \begin{eqnarray*}
&&\int_\Omega  \frac{\lambda_\beta\xi}{x_N} f\, dx = \int_\Omega \frac{\lambda_\beta\xi}{x_N} \mathcal{L}_{\beta,\varepsilon}u_\varepsilon\, dx
 \\&=& \int_\Omega  u_\varepsilon (-\Delta \frac{\lambda_\beta\xi}{x_N})\, dx+\int_{\partial\Omega}  g \frac{\partial (|x|^{\tau_+(\beta)}\xi)}{\partial \nu}\, d\omega (x)+\int_\Omega\frac{\beta}{|x|^2+\varepsilon}u_\varepsilon\lambda_\beta\xi\, dx \\
    &=& \int_\Omega  u_\varepsilon  \mathcal{L}_\beta^* (\frac{\xi}{x_N})\, d\gamma_\beta+\int_{\partial\Omega}  g \frac{\partial \xi}{\partial \nu}\, d\omega_\beta - \beta\varepsilon \int_\Omega\frac{  u_\varepsilon \xi}{(|x|^2+\varepsilon)|x|^2 x_N} \, d\gamma_\beta.
 \end{eqnarray*}

Note that if  $f\ge0$ in $\Omega$ and $g\ge0$ on $\partial\Omega$, then $u_\varepsilon\ge0$ in $\Omega$. Let $\varepsilon_1\ge\varepsilon_2$ and $u_{\varepsilon_1},\, u_{\varepsilon_2}$ be two  solutions of  (\ref{eq2.2})
respectively.
If $\beta\ge \beta_0$,  we observe that $\mathcal{L}_{\beta,\varepsilon_2}u_{\varepsilon_1}\ge \mathcal{L}_{\beta,\varepsilon_1}u_{\varepsilon_1} =f,$
so $u_{\varepsilon_1}$ is a super solution of (\ref{eq2.2}) with $\varepsilon=\varepsilon_2$ and by comparison principle, it holds
$u_{\varepsilon_1}\ge u_{\varepsilon_2}\;\; {\rm in}\; \Omega.$
The proof ends.\hfill$\Box$\smallskip

Now we  build the distributional identity for the classical solution of nonhomogeneous problem with $g=0$ and moderate singularity at the origin, i.e.
\begin{equation}\label{3.1}
 \lim_{r\to0^+} \sup_{x\in \partial_+ B_r(0)} \frac{|u(x)|}{\Lambda_\beta (x)} =0.
\end{equation}

\begin{proposition}\label{pr 2.1}
Let $\beta\ge\beta_0$, $N\ge 2$, $f\in C_{loc}^{\theta}( \bar\Omega )$ with $\theta \in(0,1)$,    then
\begin{equation}\label{2.02}
 \arraycolsep=1pt\left\{
\begin{array}{lll}
 \displaystyle \quad \mathcal{L}_\beta u= f\quad
   {\rm in}\;\;  {\Omega},\\[2mm]
 \phantom{  L_\beta   }
\quad \displaystyle  u= 0\quad  {\rm   on}\;\; \partial{\Omega}\setminus \{0\},
 \end{array}\right.
\end{equation}
subjecting to  (\ref{3.1}),
has a unique solution $u_\beta$, which 
satisfies the distributional identity
\begin{equation}\label{2.01}
 \int_{\Omega} u_\beta  \mathcal{L}_\beta^*(\frac{\xi}{x_N} )\, d\gamma_\beta = \int_{\Omega}  \frac{f\xi}{x_N} \, d\gamma_\beta,\quad \forall\,\xi\in  C^{1.1}_0(\Omega).
\end{equation}
\end{proposition}

 \noindent {\bf Proof.} The uniqueness follows by Lemma \ref{lm cp}. Since $\mathcal{L}_\beta$ is a linear operator, we only have to deal with
the case that $f\ge0$ in $\Omega$.

\vskip2mm
\noindent{\it Part 1:} $\beta>0$.  In this case, the mapping  $\varepsilon \mapsto u_\varepsilon$ is decreasing, where $u_\varepsilon>0$ is the solution of (\ref{eq2.2}) with $g=0$. Then 
$u_\beta :=\lim_{\varepsilon\to0^+} u_\varepsilon$ exists, and by the standard regularity theory,
we have that
$u_\beta$ is a classical solution of
\begin{equation}\label{2.020}
 \arraycolsep=1pt\left\{
\begin{array}{lll}
 \displaystyle   \mathcal{L}_\beta u= f\quad
   {\rm in}\;\;  {\Omega},\\[2mm]
 \phantom{  L_\beta   }
 \displaystyle  u= 0\quad  {\rm   on}\;\; \partial{\Omega}.
 \end{array}\right.
\end{equation}

\noindent{\it Part 2:} $\beta\in[\beta_0,0)$. Without loss of generality, we assume that $\Omega\subset B_\frac12(0)$.  Denote
$$
V_{t,s}(x):=\left\{
\begin{array}{lll}
 \displaystyle  t x_N|x|^{-\frac N2}-sx_N^2|x|^{\tau_+(\beta)}  \quad
   &{\rm if}\;\; \beta\in(\beta_0,0),\\[2mm]
 \phantom{     }
 \displaystyle  tx_N|x|^{-\frac N2} ( -\ln  |x|)^{\frac12}  -sx_N^2|x|^{-\frac N2}    \quad  &{\rm   if}\;\; \beta=\beta_0,
 \end{array}\right.
$$
where the parameters $s,t\ge0$.

Then for $\beta\in(\beta_0,0)$, we see that   $V_{t,s}(x)>0$ for $x\in\Omega$ if $t\ge s$ and
$$\mathcal{L}_\beta V_{t,s}(x)=tc_\beta(-  N/2)x_N|x|^{-\frac N2-2}+ 2s|x|^{\tau_+(\beta)} +2s\tau_+(\beta) x_N^2|x|^{\tau_+(\beta)-2},$$
where  $c_\beta(-  N/2)>0$ and $\tau_+(\beta)<0$.
Since $f$   is bounded in $\Omega$,   let
$$s_0=\frac12\sup_{x\in\Omega} \frac{|f(x)|}{|x|^{\tau_+(\beta)}}$$
and then we fix $t_0\ge s_0$ such that
$$t_0c_\beta(-  N/2)x_N|x|^{-\frac N2-2}+2s_0\tau_+(\beta) x_N^2|x|^{\tau_+(\beta)-2}\ge0. $$
So $V_{t_0,s_0}$ is a positive supersolution of (\ref{2.02}).

For $\beta=\beta_0, \, \tau_-(\beta) = - \frac{N}{2}$, we have that
$$\mathcal{L}_\beta V_{t,s}(x)= \frac t4 x_N|x|^{-\frac N2-2} (-\ln |x|)^{ -\frac12}  + 2s|x|^{-\frac N2} -2sN x_N^2|x|^{-\frac N2-2}.$$
We take $s_0$ as above where $\beta$ is replaced by $\beta_0$
and we fix $t_0\ge s_0$ such that
$$\frac{t_0}{4} x_N|x|^{-\frac N2-2} (-\ln |x|)^{ -\frac12}  -2s_0N x_N^2|x|^{-\frac N2-2}\ge0. $$
So $V_{t_0,s_0}$ is also a positive supersolution of (\ref{2.02}) in this case which implies,
by comparison principle, that we have
$$u_{\varepsilon}(x)\le V_{t_0,s_0}(x),\quad\forall\, x\in\Omega.$$

\noindent{\it Proof of (\ref{2.01}).} We need to estimate $ \displaystyle\int_\Omega\frac{  u_\varepsilon \xi}{(|x|^2+\varepsilon)|x|^2 x_N}  \, d\gamma_\beta $ for $0 < \varepsilon < \varepsilon_0$ for some $\varepsilon_0>0$ fixed.  we first consider the case $\beta>0$.
We observe that
\begin{eqnarray*}
&&\varepsilon\int_{\Omega\setminus B_{  \sqrt{\varepsilon} }(0)} \frac{ u_\varepsilon \xi  \lambda_\beta(x)  }{(|x|^2+\varepsilon)|x|^2x_N}\, dx
\\ &\le &\varepsilon\norm{u_{\varepsilon_0}}_{L^\infty(\Omega)}\norm{\xi/\rho}_{L^\infty(\Omega)}\int_{\Omega\setminus B_{ \sqrt{\varepsilon}}(0)} \frac{ |x|^{\tau_+(\beta)-2} } { |x|^2+\varepsilon }   \, dx   \\
  &\le & \norm{u_{\varepsilon_0}}_{L^\infty(\Omega)}\norm{\xi/\rho}_{L^\infty(\Omega)} \varepsilon^{\frac{N-2+\tau_+(\beta)}2} \int_{B_{\frac{1}{2\sqrt{\varepsilon}}}(0)\setminus B_1(0)}   |y|^{\tau_+(\beta)-4}   \, dy
  \\&\le& c_3\norm{u_{\varepsilon_0}}_{L^\infty(\Omega)}\norm{\xi/\rho}_{L^\infty(\Omega)}(2^{-\tau_+(\beta)+4-N}\varepsilon+\varepsilon^{\frac{N-2+\tau_+(\beta)}2})
  \\[2mm] &\to&0\ \quad{\rm as}\quad \ \varepsilon\to0^+
\end{eqnarray*}
and
\begin{eqnarray*}
&&\varepsilon\int_{ B_{ \sqrt{\varepsilon }}(0)} \frac{ u_\varepsilon \xi \lambda_\beta(x)  }{(|x|^2+\varepsilon)|x|^2x_N}\,   \, dx
\\ &\le &  \norm{u_{\varepsilon_0}}_{L^\infty(\Omega)}\norm{\xi/\rho}_{L^\infty(\Omega)}\int_{  B_{ \sqrt{\varepsilon }}(0)}   |x|^{\tau_+(\beta)-2}      \, dx   \\
  &\le & c_4\norm{u_{\varepsilon_0}}_{L^\infty(\Omega)}\norm{\xi/\rho}_{L^\infty(\Omega)} \varepsilon^{\frac{N-2+\tau_+(\beta)}2}
   \\[2mm] &\to&0\ \quad{\rm as}\quad\ \varepsilon\to0^+,
\end{eqnarray*}
 where $\rho(x)={\rm dist}(x,\,\partial\Omega)$ and $\frac{N-2+\tau_+(\beta)}2>0$.
Therefore, passing to the limit of (\ref{11.1}), we obtain (\ref{2.01}).

For $\beta\in(\beta_0,\,0)$,   from the increasing monotonicity and the upper bound $V_{s_0,t_0}$, we have that
$$\lim_{\varepsilon\to0^+}\int_\Omega  u_\varepsilon  \mathcal{L}_\beta^* (\frac{\xi}{x_N}) d\gamma_\beta=\int_\Omega  u_\beta \mathcal{L}_\beta^* (\frac{\xi}{x_N}) d\gamma_\beta $$
and
$$\varepsilon  \int_\Omega\frac{ \xi u_\varepsilon\lambda_\beta(x)}{(|x|^2+\varepsilon)|x|^2x_N}  dx\le c_5\varepsilon \int_\Omega \frac{ |x|^{-N+ \sqrt{\beta-\beta_0}}}{|x|^2+\varepsilon}  dx.  $$
By directly compute, we have that
\begin{eqnarray*}
\varepsilon \int_{\Omega\setminus B_{ \sqrt{\varepsilon}}(0)} \frac{ |x|^{-N+ \sqrt{\beta-\beta_0}}}{|x|^2+\varepsilon}  dx
  &\le & c_6 \varepsilon^{\frac{ \sqrt{\beta-\beta_0}}2} \int_{B_{\frac{1}{2\sqrt{\varepsilon}}}(0)\setminus B_1(0)}   |y|^{-N-2+\sqrt{\beta-\beta_0}}   \, dy
  \\&\le& c_7 ( \varepsilon+\varepsilon^{\frac{ \sqrt{\beta-\beta_0}}2} ) \to0\quad{\rm as}\quad \varepsilon\to0^+
\end{eqnarray*}
and
\begin{eqnarray*}
\varepsilon\int_{ B_{ \sqrt{\varepsilon }}(0)} \frac{ |x|^{-N+ \sqrt{\beta-\beta_0}}}{|x|^2+\varepsilon}  dx   &\le &  \int_{  B_{ \sqrt{\varepsilon }}(0)}   |x|^{-N+ \sqrt{\beta-\beta_0}}     \, dx   \\
  &\le & c_8  \varepsilon^{\frac{ \sqrt{\beta-\beta_0}}2} \to 0\quad{\rm as}\quad \varepsilon\to0^+,
\end{eqnarray*}

As a conclusion, passing to the limit in (\ref{11.1}) as $\varepsilon\to0^+$, we have that $u_\beta$ satisfies that
\begin{equation}\label{2.01-0}
 \int_{\Omega} u_\beta  \mathcal{L}_\beta^*(\frac{\xi}{x_N} )\, d\gamma_\beta = \int_{\Omega}   \frac{f\xi}{x_N} \, d\gamma_\beta,\quad \forall\,\xi\in  C^{1.1}_0(\Omega).
\end{equation}

\noindent Finally, we prove (\ref{2.01}) with $\beta=\beta_0$, We  claim that
the mapping $\beta\mapsto u_\beta$ with $\beta\in(\beta_0,0)$ is decreasing.
In fact, if $\beta_0<\beta_1\le\beta_2<0$, we know that
\begin{eqnarray*}
f=\mathcal{L}_{\beta_1} u_{\beta_1} &=& -\Delta u_{\beta_1}+\frac{\beta_1}{|x|^2}  u_{\beta_1} \\
    &\le& -\Delta u_{\beta_1}+\frac{\beta_2}{|x|^2}  u_{\beta_1}= \mathcal{L}_{\beta_2} u_{\beta_1},
\end{eqnarray*}
by Lemma \ref{lm cp},
which implies that $u_{\beta_1}\ge u_{\beta_2}.$

We know that $V_{s_0,t_0}$ is a super solution of (\ref{2.02})  with $\beta\in(\beta_0,0)$. So it follows by Lemma \ref{lm cp} that
$\{u_\beta\}_\beta$ is uniformly bounded by the upper bound $V_{s_0,t_0} \in L^1(\Omega,\frac1{x_N}d\gamma_\beta)$.

For $\xi \in  C^{1.1}_0(\Omega)$, we have  that
$$|\mathcal{L}_\beta^*(\frac{\xi}{x_N})|\le c_{9}(\|\frac{\xi}{x_N}\|_{C^{1.1}(\Omega)}+  \|\frac{\xi}{x_N}\|_{C^{1}(\Omega)} x_N^{-1}), $$
where $c_{9}>0$ is independent of $\beta$.

From the dominate monotonicity convergence theorem and the uniqueness of the solution, we have that
$$u_{\beta} \to u_{\beta_0} \quad{\rm a.e.\ in}\  \Omega \; {\rm as}\;\; \beta\to \beta^+_0 \quad {\rm and \; in}\;\; L^1(\Omega,\, x_N^{-1}d\gamma_\beta)$$
and $u_{\beta_0}$ is a classical solution of (\ref{2.02}) with $\beta=\beta_0$.
 Passing to the limit of (\ref{2.01-0}) as $\beta\to \beta^+_0$ to obtain that
$$  \int_{\Omega}u_{\beta_0}  \mathcal{L}^*_{\beta_0}(\frac{\xi}{x_N})\,  d\gamma_{\beta_0}= \int_{\Omega}\frac{f \xi}{x_N}  d\gamma_{\beta_0}.$$
The proof ends.\hfill$\Box$

\begin{remark}

We note that when $\beta\ge 0$ and $f$ is bounded, the moderate singular solution of problem (\ref{2.02})  is no longer singular, that means, it is a classical solution of
\begin{equation}\label{2.0200}
 \arraycolsep=1pt\left\{
\begin{array}{lll}
 \displaystyle  \mathcal{L}_\beta u= f\quad
   {\rm in}\;\;  {\Omega},\\[1.5mm]
 \phantom{  L_\beta  }
 \displaystyle  u= 0\quad  {\rm   on}\;\; \partial{\Omega}.
 \end{array}\right.
\end{equation}
\end{remark}

Now we prove the following
\begin{lemma}\label{lm 2.5}
$(i)$ The problem
\begin{equation}\label{eq 2.2}
\arraycolsep=1pt\left\{
\begin{array}{lll}
 \displaystyle   \mathcal{L}_{\beta}^* (\frac{u}{x_N}) = 1\quad
   {\rm in}\;\;   \Omega,\\[2mm]
 \phantom{ \mathcal{L}_{\beta}^*-\ \,\,\, }
 \displaystyle  u= 0\quad  {\rm   on}\;\; \partial{\Omega}
 \end{array}\right.
\end{equation}
has a unique positive solution $w_1\in C^2(\Omega)\cap C^{0.1}_0(\Omega)$.

$(ii)$ The problem
\begin{equation}\label{eq 2.2-1}
\arraycolsep=1pt\left\{
\begin{array}{lll}
 \displaystyle   \mathcal{L}_{\beta}^* (\frac{u}{x_N}) = \frac{1}{x_N} \quad
   &{\rm in}\;\;   \Omega,\\[2mm]
 \phantom{ \mathcal{L}_{\beta}^*--  }
 \displaystyle  u= 0\quad  &{\rm   on}\;\; \partial{\Omega}
 \end{array}\right.
\end{equation}
has a unique positive solution $w_2\in C^2(\Omega)\cap C^{1}_0(\bar\Omega\setminus\{0\})\cap C^{0.1}_0(\Omega)$.

\end{lemma}
 \noindent {\bf Proof.}
We first claim that problem (\ref{2.02}) has a unique classical positive solution $w_\beta$ under the constraint (\ref{3.1})
when $f(x)=\lambda_\beta(x)$ or $f(x)=|x|^{\tau_+(\beta)}$.

In fact, let $f_n(x)=\lambda_\beta(x)\eta_0(n|x|)$, where $\eta_0:[0,+\infty)\to [0,\,1]$ is  a decreasing $C^\infty$ function
 satisfying (\ref{eta}). Then $f_n\in C^{\theta}( \bar\Omega )$ with $\theta \in(0,1), f_n \leq f$,  and by Proposition \ref{pr 2.1}, let $w_n$ be the solution of
 problem
 \begin{equation}\label{2.02-00}
 \arraycolsep=1pt\left\{
\begin{array}{lll}
 \displaystyle  \mathcal{L}_\beta u= f_n\quad
   &{\rm in}\;\;  {\Omega},\\[2mm]
 \phantom{  L_\beta   }
 \displaystyle  u= 0\quad  &{\rm   on}\;\; \partial{\Omega}\setminus \{0\},
 \end{array}\right.
\end{equation}
subject to  (\ref{3.1}). We know that the mapping: $n\to w_n$ is increasing by the increasing monotone of $\{f_n\}$.
So we only construct a suitable upper bound for $w_n$ in the cases that $f(x)=\lambda_\beta(x)$ and $f(x)=|x|^{\tau_+(\beta)}$
respectively.

\vskip2mm
When $f(x)=\lambda_\beta(x)$, let $ V_{t,s}(x)= t \lambda_\beta(x)-sx_N|x|^{\tau_+(\beta)+2} $
for $s,t>0$.
It is know that
 $$ \mathcal{L}_\beta V_{t,s}=-sc_{\tau_+(\beta)+2}\lambda_\beta(x),\quad x\in \R^N_+, $$
for some $c_{\tau_+(\beta)+2}<0$.  So fix $s= -1/c_{\tau_+(\beta)+2}$ and then fix $t>0$ such that
 $$V_{t,s}(x)>0,\quad \forall x\in\Omega.$$
The limit of $\{w_n\}_n$, denoting by $w_{\beta,1}$, is a solution of (\ref{3.1})
satisfying
$w_{\beta,1}\le  V_{t,s}(x).$

\vskip2mm
 When $f(x)=|x|^{\tau_+(\beta)}$, let
$$ W_{t,s,l}(x)= t \lambda_\beta(x)-s(x_N|x|^{\tau_+(\beta)+2}+lx_N^2|x|^{\tau_+(\beta)+2}), $$
where $s,t,l>0$.
 We observe that
 $$ \mathcal{L}_\beta W_{t,s,l}(x)= s[-c_{\tau_+(\beta)+2}  \lambda_\beta(x) +  2l|x|^{\tau_+(\beta)} +2l\tau_+(\beta) x_N^2|x|^{\tau_+(\beta)}],\; x\in \R^N_+, $$
 with the same constant $c_{\tau_+(\beta)+2}<0$ as above.  Then we choose $l>0$ such that $-2c_{\tau_+(\beta)+2}l\tau_+(\beta) x_N>0$ for $x\in\Omega$, $s=\frac1{2l}$ and we take $t>0$ such that
 $W_{t,s,l}>0$ in $\Omega$ and
 $$\mathcal{L}_\beta W_{t,s,l}(x)\ge  |x|^{\tau_+(\beta)}.$$
Thus, the limit of $\{w_n\}_n$, denoting by $w_{\beta,2}$, is a solution of (\ref{3.1})
such that
$$w_{\beta,2}(x)\le  W_{t,s,l}(x).$$
As a conclusion, for $i=1,2$,
\begin{equation}\label{0.3}
 w_{\beta,i}\le t\lambda_\beta\quad{\rm in}\quad \Omega.
\end{equation}

Denote $w_i= w_{\beta,i}x_N/\lambda_\beta,$
we observe that
\begin{eqnarray*}
1 =  \lambda_\beta^{-1}  \mathcal{L}_\beta w_{\beta,1} = \lambda_\beta^{-1}  \mathcal{L}_\beta (\lambda_\beta w_{1}/x_N)
   = \mathcal{L}_\beta^* (w_1/x_N)
\end{eqnarray*}
and
\begin{eqnarray*}
1/x_N =  \lambda_\beta^{-1}  \mathcal{L}_\beta w_{\beta,2} = \lambda_\beta^{-1}  \mathcal{L}_\beta (\lambda_\beta w_{2}/x_N)
   = \mathcal{L}_\beta^* (w_2/x_N).
\end{eqnarray*}
Moreover, by (\ref{0.3}), it follow that $w_i\le t x_N.$
Then we have that $w_i\in C^2(\Omega)\cap C^{0.1}_0(\Omega)$ for $i=1,2$.  Away from the origin,
Hardy's operator is uniform elliptic,  thus $u \in C^1_0(\bar\Omega\setminus\{0\})$ and then
$u\in C^2(\Omega)\cap C^{1}_0(\bar\Omega\setminus\{0\})\cap C^{0.1}_0(\Omega).$
\hfill$\Box$

\medskip
Although $C^2(\Omega)\cap C^{1}_0(\bar\Omega\setminus\{0\})\cap C^{0.1}_0(\Omega)$ is not suitable as test function space for problem (\ref{eq 1.1fg}),  $w_1,\, w_2$ are
still valid as test functions for formula (\ref{1.2fg}) with $k=0$ in the distributional sense.

\vskip2mm
For given $f\in C^1(\bar\Omega)$, a direct consequence of Lemma \ref{lm 2.5} can be stated as follows
\begin{corollary}\label{re 3.1}
Assume that  $f\in C^1(\bar\Omega\setminus\{0\})$ satisfying for some $c_{10}>0$
$$|f(x)|\le \frac {c_{10}}{x_N}.$$
Then there exists a unique solution of $w_f\in C^2(\Omega)\cap C^{0.1}_0(\Omega)$ of
\begin{equation}\label{eq 2.2f}
\arraycolsep=1pt\left\{
\begin{array}{lll}
 \displaystyle   \mathcal{L}_{\beta}^* (\frac{u}{x_N}) = f\quad
  &{\rm in}\;\;   \Omega,\\[2mm]
 \phantom{ \mathcal{L}_{\beta}^*-\ \,\,\, }
 \displaystyle  u = 0\quad  &{\rm   on}\;\; \partial{\Omega}.
 \end{array}\right.
\end{equation}

\end{corollary}

\setcounter{equation}{0}
\section{Fundamental solution}

\subsection{ In half space}

In this subsection, we give the proof of Theorem \ref{teo 0}.
\smallskip

 \noindent {\bf Proof of Theorem \ref{teo 0}.}
 For any $\xi\in C^{1.1}_0(\R^N_+)$,
 we know there exists a unique $\zeta\in C^{1.1}_c(\ \R^N)$ such that
  $\xi(x)=x_N\zeta(x)$ for $x\in\overline{\R^N_+}$. Moreover, we have that $\frac{\partial \xi}{\partial x_N}(0)=\zeta(0)$.

  Take $\zeta\in  C^{1.1}_c( \R^N)$, multiplying $\lambda_\beta\zeta$ in (\ref{eq 1.1}) and integrating over $\R^N_+\setminus \overline{B_r(0)}$, then we have that
\begin{eqnarray*}
0 &=& \int_{\R^N_+\setminus \overline{B_r(0)}} \mathcal{L}_\beta (\Lambda_\beta) \lambda_\beta \zeta\, dx
    =   \int_{\R^N_+\setminus \overline{B_r(0)}}\Lambda_\beta \mathcal{L}_\beta^*(\zeta )\, d\gamma_\beta  \\&&
   +\int_{\partial_+ B_r(0)}\Big(- \nabla \Lambda_\beta\cdot\frac{x}{|x|}  \lambda_\beta  +\nabla  \lambda_\beta  \cdot\frac{x}{|x|}\Lambda_\beta\Big)\zeta\,d\omega
 \\&& + \int_{\partial_+ B_r(0)} \Lambda_\beta \lambda_\beta  \Big( \nabla\zeta\cdot\frac{x}{|x|}\Big) \,d\omega,
\end{eqnarray*}
where $\partial_+ B_r(0)=\partial B_r(0)\cap \R^N_+$.
For $\beta \geq \beta_0$,  we see that for $r=|x|>0$ small,
\begin{eqnarray*}
&& -\nabla \Lambda_\beta(x)\cdot\frac{x}{|x|} \lambda_\beta (x) +\nabla \lambda_\beta(x) \cdot\frac{x}{|x|}\Lambda_\beta(x)\\[2mm]&=&\left\{\arraycolsep=1pt
\begin{array}{lll}
   2\sqrt{\beta-\beta_0} \,x_N^2r^{-N-1}\quad
   &{\rm if}\;\; \beta>\beta_0,\\[1.5mm]
 \phantom{   }
 x_N^2r^{-N-1} \quad  &{\rm  if}\;\; \beta=\beta_0
 \end{array}
 \right.
\end{eqnarray*}
and
$$ |\zeta(x) - \zeta(0)|\le c_{11}r,  $$
then
\begin{eqnarray*}
 && \int_{\partial_+ B_r(0)} \sqrt{\beta-\beta_0} \, x_N^2r^{-N-1}   \zeta(0) x_N d\omega(x)
  \\[2mm]& =& \left\{\arraycolsep=1pt
\begin{array}{lll}
   \sqrt{\beta-\beta_0}\, \displaystyle\int_{\partial_+ B_1(0)} x_N^2 d\omega(x)\,\zeta(0) \quad
   &{\rm if}\;\; \beta>\beta_0,\\[2mm]
 \phantom{   }
\displaystyle\int_{\partial_+ B_1(0)} x_N^2 d\omega(x) \, \zeta(0)\quad  &{\rm  if}\;\; \beta=\beta_0
 \end{array}
 \right.
\\[2mm]  & = & c_\beta \zeta(0)
\end{eqnarray*}
and
\begin{eqnarray*}
 &&\Big| \int_{\partial_+ B_r(0)} \Big(- \nabla \Lambda_\beta\cdot\frac{x}{|x|}  \lambda_\beta  +\nabla  \lambda_\beta  \cdot\frac{x}{|x|}\Lambda_\beta\Big)\zeta\,d\omega
 -c_\beta \zeta(0)  \Big|
 \\[2mm]&\le& c_{12}(\sqrt{\beta-\beta_0} +1)\,r\int_{\partial_+ B_1(0)} x_N^2 d\omega(x)
 \\[1.5mm]&\to& 0\quad{\rm as}\quad r\to0^+,
\end{eqnarray*}
that is,
$$
 \lim_{r\to0}\Big(\int_{\partial_+ B_r(0)} - \nabla \Lambda_\beta\cdot\frac{x}{|x|}  \lambda_\beta\zeta \,d\omega
   + \int_{\partial B_r(0)} \nabla  \lambda_\beta \cdot\frac{x}{|x|}\Lambda_\beta\zeta \,d\omega\Big)= c_\beta \zeta(0).
$$
Moreover, we see that
$$
 \Big|\int_{\partial_+ B_r(0)} \Lambda_\beta \lambda_\beta  \Big( \nabla\zeta\cdot\frac{x}{|x|}\Big) \,d\omega \Big| \le  \norm{\zeta}_{C^1}\, r\int_{\partial_+ B_1(0)} x_N^2  d\omega
\to 0\quad{\rm as}\quad r\to0^+. $$
Therefore, we have that
$$\lim_{r\to0^+}\int_{\R^N\setminus \overline{B_r(0)}}\Lambda_\beta \mathcal{L}_\beta^*( \zeta ) d\gamma_\beta=c_\beta   \zeta(0),$$
which implies (\ref{1.2}). The proof ends. \hfill$\Box$

\subsection{ Trace of $\Lambda_\beta$.}
The following theorem shows the trace of $\Lambda_\beta$.
\begin{theorem}\label{teo 2.2}
Let $d\omega_\beta(x') = |x'| ^{\tau_+(\beta)} dx'$
for $x'\in\R^{N-1}$,
then for any $\zeta \in C_c(\R^{N-1})$,
\begin{equation}\label{2.2}
\lim_{t\to0^+} \int_{\R^{N-1}}\Lambda_\beta(x',t)   \zeta (x')   d\omega_\beta(x')  =b_N  \zeta (0),
\end{equation}
where    $$b_N=\int_{\R^{N-1}}  (1+|y'|^2)^{-\frac N2} dy'>0.$$
 This is to say that the trace of $\Lambda_\beta$ is $\delta_0$ in the
$d\gamma_\beta$-distributional sense.

 \end{theorem}
\noindent {\bf Proof.} For any $\zeta \in C_c(\R^{N-1})$, there exists $R>0$ such that supp$\,\zeta\subset B_R'(0)$, here and in the sequel, denoting by $B_R'(0)$ the ball in $\R^{N-1}$.
By direct computations, we have that
\begin{eqnarray*}
 \int_{\R^{N-1}}\Lambda_\beta(x',t)   \zeta (x')\,  d\omega_\beta(x') & =&\int_{B_R'(0)}\Lambda_\beta(x', t)   \zeta (x')  \, d\omega_\beta(x')
 \\&=&  \int_{B_{R/t}'(0)} (|y'|^2+1)^{\frac{\tau_-(\beta)}2}|y'|^{\tau_+(\beta)} \zeta (ty')  dy'.
\end{eqnarray*}
For any $\varepsilon>0$, there exists $R_\varepsilon>1$ such that
 \begin{eqnarray*}
 &&\int_{B_{R/t}'(0)\setminus B'_{R_\varepsilon }(0)} (|y'|^2+1)^{\frac{\tau_-(\beta)}2}|y'|^{\tau_+(\beta)} \zeta (ty')  dy'
 \\&\le& \norm{\zeta}_{L^\infty(\R^{N-1})}\int_{\R^{N-1}\setminus B'_{R_\varepsilon }(0)}  |y'| ^{-N}  dy'\\&\le & \norm{\zeta}_{L^\infty(\R^{N-1})} |\mathcal{S}^{N-2}|\varepsilon,
 \end{eqnarray*}
 where $R_\varepsilon\le \frac1\varepsilon$. Let
 $$A:=\int_{  B'_{R_\varepsilon}(0)} (|y'|^2+1)^{\frac{\tau_-(\beta)}2}|y'|^{\tau_+(\beta)} \zeta (ty')  dy'-\int_{\R^{N-1}} (|y'|^2+1)^{\frac{\tau_-(\beta)}2}|y'|^{\tau_+(\beta)} \zeta(0) dy',$$
we have that
 \begin{eqnarray*}
 | A| &\le&   \int_{  B'_{R_\varepsilon}(0)} (|y'|^2+1)^{\frac{\tau_-(\beta)}2}|y'|^{\tau_+(\beta)}\left|\zeta (ty')-\zeta(0)\right|\,dy'
 +\varepsilon |\zeta(0)|  |\mathcal{S}^{N-2}|
 \\&\le &  t \norm{\zeta}_{C^1(\R^{N-1})}\int_{  B'_{R_\varepsilon}(0)} (|y'|^2+1)^{\frac{\tau_-(\beta)}2}|y'|^{\tau_+(\beta)} dy'+\varepsilon |\zeta(0)| |\mathcal{S}^{N-2}|
  \\&= & R_\varepsilon t\norm{\zeta}_{C^1(\R^{N-1})}+\varepsilon |\zeta(0)| |\mathcal{S}^{N-2}|
  \\[1.1mm]&\le &\left(\norm{\zeta}_{C^1(\R^{N-1})}+|\zeta(0)| |\mathcal{S}^{N-2}|\right)\varepsilon,
 \end{eqnarray*}
if we take $t=\varepsilon^2$.
Passing to the limit as $\varepsilon\to0$,   we derive (\ref{2.2}).\hfill$\Box$

\subsection{Fundamental solution in bounded domain}

In this subsection, we do an approximation of the isolated singular solution.

\begin{proposition}\label{pr 2.2}
Let $\Omega$ be a $C^2$ domain  verifying (\ref{O}).
Then the problem
\begin{equation}\label{2.2.1}
 \arraycolsep=1pt\left\{
\begin{array}{lll}
 \displaystyle  \mathcal{L}_\beta u= 0\quad
   {\rm in}\;\;  {\Omega},\\[2mm]
 \phantom{  L_\beta  }
 \displaystyle  u= 0\quad  {\rm on}\;\; \partial{\Omega}\setminus \{0\},\\[2mm]
 \phantom{   }
  \displaystyle \lim_{r\to0^+}\sup_{x\in B_r^+(0)}\frac{|u(x)-\Lambda_\beta(x)|}{\Lambda_\beta (x)}=0
 \end{array}\right.
\end{equation}
admits a unique solution $\Lambda^\Omega_\beta$ satisfying
the following distributional identity:
\begin{equation}\label{2.01b}
 \int_{\Omega} \Lambda^\Omega_\beta  \mathcal{L}_\beta^*(\frac{\xi}{x_N} )\, d\gamma_\beta = c_\beta\frac{\partial \xi}{\partial x_N}(0),\quad \forall\,\xi\in  C^{1.1}_0(\Omega).
\end{equation}
\end{proposition}
\noindent {\bf Proof.}   Let $\eta_{r_0}(t)=\eta_0(\frac2{r_0}t)$, which satisfies  that
\begin{equation}\label{eta0}
 \eta_{r_0}=1\quad{\rm in}\quad [0,r_0/2]\quad{\rm and}\quad \eta_{r_0}=0\quad {\rm in}\quad[r_0,+\infty).
\end{equation}
 For $i=1,2$ the problem
\begin{equation}\label{2.3}
\arraycolsep=1pt\left\{
\begin{array}{lll}
 \displaystyle  \mathcal{L}_\beta w_i= -\nabla \eta_{r_0}\cdot \nabla  \Lambda_\beta-\Lambda_\beta \Delta\eta_{r_0}\quad
   {\rm in}\;\;  {\Omega},\\[2mm]
 \phantom{  L_\beta   }
 \displaystyle  w_i= 0\quad  {\rm   on}\;\; \partial{\Omega}\setminus \{0\},\\[2mm]
 \phantom{   }
  \displaystyle \lim_{e\in \mathcal{S}^N_+,\,t\to0^+}w_i(te)\Lambda_\beta^{-1}(te)= 2-i,
 \end{array}\right.
\end{equation}
admits a unique solutions $w_1$ and $w_2$ respectively. Obviously,
$$w_1=\Lambda_\beta\eta_{r_0}$$
and
  $-\nabla \eta_{r_0}\cdot \nabla  \Lambda_\beta-\Lambda_\beta \Delta\eta_{r_0}$ has compact set in $\Omega\cap (\overline{B_{r_0}(0)\setminus B_{\frac{r_0}2}(0)})$
and then  $-\nabla \eta_{r_0}\cdot \nabla  \Lambda_\beta-\Lambda_\beta \Delta\eta_{r_0}$ is smooth and bounded,  it follows by the proof of
Proposition \ref{pr 2.1}  that there exist $s_0,t_0>0$ such that
$|w_2|\le V_{s_0,t_0}$.

\vskip2mm
For $i=1$, following the proof of Theorem \ref{teo 0}, we get then for any $\xi\in C^{1.1}_0(\Omega)$,
\begin{equation}\label{2.4}
 \int_{\Omega}w_1 \mathcal{L}_\beta^*( \frac{\xi}{x_N})\, d\gamma_\beta =\int_{\Omega} \Big(-\nabla \eta_{r_0}\cdot \nabla
  \Lambda_\beta-\Lambda_\beta \Delta\eta_{r_0}\Big)\frac{\xi}{x_N}\, d\gamma_\beta+c_\beta \frac{\partial\xi}{\partial x_N}(0).
\end{equation}
For $i=2$,  it follows by Proposition \ref{pr 2.1}   that for any $\xi\in C^{1.1}_0(\Omega)$,
\begin{equation}\label{2.5}
 \int_{\Omega} w_2 \mathcal{L}_\beta^*( \frac{\xi}{x_N}) d\gamma_\beta=\int_{\Omega} \Big(-\nabla \eta_{r_0}\cdot \nabla  \Lambda_\beta-\Lambda_\beta \Delta\eta_{r_0}\Big) \frac{\xi}{x_N}\,d\gamma_\beta.
\end{equation}
Let $\Lambda^\Omega_\beta=\Lambda _\beta\eta_{r_0}-w_2,$ it follows by  (\ref{2.4}) and (\ref{2.5}) that
$$ \int_{\Omega}\Lambda^\Omega_\beta \mathcal{L}_\beta^*( \frac{\xi}{x_N})\, d\gamma_\beta =c_\beta \frac{\partial\xi}{\partial x_N}(0),\quad \forall\,\xi\in  C^{1.1}_0(\Omega).$$

Finally, it's clear that if $u_1$ and $u_2$ are two solutions of (\ref{2.2.1}), then $w := u_1 - u_2$ satisfies
$$\lim_{r\to0^+}\sup_{x\in B_r^+(0)}\frac{|w(x)|}{\Lambda_\beta (x)}=0.$$
Combining with the fact that
$$ \mathcal{L}_\beta w= 0\;\;
   {\rm in}\;\;  {\Omega}\quad{\rm and}\quad
   w= 0\;\;  {\rm on}\;\; \partial{\Omega}\setminus \{0\},$$
and Lemma \ref{lm cp}, we have that $w\equiv0$. Thus the uniqueness is proved.\hfill$\Box$

 \smallskip

\setcounter{equation}{0}
\section{Existence}

\subsection{Zero Dirichlet boundary}

Our purpose in this section is to clarify the isolated singularities of   the nonhomogeneous problem
 \begin{equation}\label{eq 1.1f}
 \arraycolsep=1pt\left\{
\begin{array}{lll}
 \displaystyle   \mathcal{L}_\beta u= f\quad
   {\rm in}\;\;  \Omega,\\[1.5mm]
 \phantom{   L_\beta   }
 \displaystyle  u=0\quad  {\rm   on}\;\; \partial{\Omega}\setminus \{0\},
 \end{array}\right.
\end{equation}
where $f\in C^\theta_{loc}(\bar\Omega\setminus\{0\})$ with $\theta \in(0,1)$. Recall that $\mathcal{L}^*_\beta$  is given by (\ref{L})
and $d\gamma_\beta(x) =\lambda_\beta(x) dx$. We prove the following

\begin{theorem}\label{teo 3.1}
$(i)$ Assume  that $f\in L^1(\Omega,\, d\gamma_\beta)$ and $u\in L^1(\Omega,\frac1{|x|}d\gamma_\beta)$ is a  classical solution of problem (\ref{eq 1.1f}), then
 there exists some $k\in \R$ such that there holds
 \begin{equation}\label{1.2f}
 \int_{\Omega}u\,  \mathcal{L}_\beta^*(\frac{\xi}{x_N})\, d\gamma_\beta  = \int_{\Omega}   \frac{f\xi}{x_N} \, d\gamma_\beta +  k \frac{\partial \xi}{\partial x_N}(0) ,\quad\forall\, \xi\in   C^{1.1}_0(\Omega).
\end{equation}

$(ii)$ Inversely, assume  that $f\in L^1(\Omega,\, d\gamma_\beta)$, then  for any $k\in\R$,  problem (\ref{eq 1.1f})
has a unique  solution $u_k\in L^1(\Omega,\frac1{|x|}d\gamma_\beta)$  verifying (\ref{1.2f}) with such $k.$

\end{theorem}
\noindent{\bf Proof.} $(i)$ Let $\tilde \Omega$ be the interior set of $\bar \Omega\cup \overline{\{(x',-x_N):(x',x_N)\in\Omega\}}$ and extend $u$ (resp. $f$) by the $x_N$-odd extension to $\tilde u$ (resp. $\tilde f$) in $\tilde \Omega$, then $\mathcal{L}_\beta \tilde u=\tilde f$.
Our aim is to see the distributional property at the origin.
Denote  by $L$  the operator related to $\mathcal{L}_\beta \tilde u -\tilde f$ in the  distribution sense,  i.e.
\begin{equation}\label{3.1.1}
L(\zeta)=\int_{\tilde\Omega} \Big(\tilde u \mathcal{L}_\beta^*(\zeta) -   \tilde f\zeta \Big)|x_N||x|^{\tau_+(\beta)}\,dx,\quad \forall\zeta\in C^\infty_c(\tilde\Omega).
\end{equation}

For any $\zeta \in C^\infty_c(\tilde\Omega\setminus\{0\} )$, we have that $L(\zeta)=0.$
In fact, there exists $\varepsilon>0$ such that supp$(\zeta)\subset \tilde\Omega\setminus B_\varepsilon(0)$ and then
\begin{eqnarray*}
0&=&2\int_{  \Omega}  \zeta(\mathcal{L}_\beta  u-f)\,d\gamma_\beta
=\int_{\tilde \Omega}  \zeta(\mathcal{L}_\beta \tilde u-\tilde f) \,d\tilde \gamma_\beta
\\&=& -\int_{\tilde \Omega}  \tilde f \zeta\, d\tilde\gamma_\beta+ \int_{\Omega\setminus B_\varepsilon(0)} u \mathcal{L}_\beta^* \zeta d  \gamma_\beta+  \int_{\partial(\Omega\setminus B_\varepsilon(0))\cap(\R^{N-1}\times\{0\}) }\frac{\partial u}{\partial x_N} \zeta d\omega_\beta \\
  &&+\int_{(-\Omega)\setminus B_\varepsilon(0)} (-u) \mathcal{L}_\beta^* \zeta d\tilde\gamma_\beta  +  \int_{\partial(-\Omega\setminus B_\varepsilon(0))\cap(\R^{N-1}\times\{0\}) }\frac{\partial \tilde u}{\partial(- x_N)} \zeta d\omega_\beta   \\
  &=&\int_{\tilde \Omega\setminus B_\varepsilon(0)}  (\tilde u \mathcal{L}_\beta^*  \zeta-\tilde f\zeta)\, d\tilde\gamma_\beta
  \\&=&\int_{\tilde \Omega}  (\tilde u \mathcal{L}_\beta^*  \zeta-\tilde f\zeta)\, d\tilde\gamma_\beta  ,
\end{eqnarray*}
where $d\tilde \gamma_\beta  = |\tilde \lambda_\beta(x)|dx$, $\tilde \lambda_\beta$ is the odd extension of $\lambda_\beta$ and
 $$\int_{\partial(\Omega\setminus B_\varepsilon(0))\cap(\R^{N-1}\times\{0\}) }\frac{\partial u}{\partial x_N} \zeta d\omega_\beta=- \int_{\partial(-\Omega\setminus B_\varepsilon(0))\cap(\R^{N-1}\times\{0\}) }\frac{\partial \tilde u}{\partial(- x_N)} \zeta d\omega_\beta.$$

 By Theorem XXXV in \cite{S} (see also Theorem 6.25 in \cite{R}), it implies that
\begin{equation}\label{S}
 L=\sum_{|a|=0}^p k_a D^{a}\delta_0,
\end{equation}
where $p\in\N$, $a=(a_1,\cdots,a_N)$ is a multiple index with $a_i\in\N$, $|a|= \sum_{i=1}^Na_i$ and in particular, $D^0\delta_0=\delta_0$.
Then we have that
\begin{equation}\label{3.3}
L(\zeta)= \int_{\tilde \Omega} \Big(\tilde u \mathcal{L}_\beta^*\zeta -f\zeta\Big)\,d\tilde\gamma_\beta=\sum_{|a|=0}^{\infty} k_a D^{a}\zeta(0),\quad\ \ \forall \zeta\in C^\infty_c (\tilde \Omega).
\end{equation}

For any multiple index $a=(a_1,\cdots,a_N)$, let $\zeta_a$ be a    $C^\infty$ function such that
\begin{equation}\label{3.2-2}
{\rm supp}(\zeta_a)\subset \overline{B_2(0)}\quad{\rm and}\quad \zeta_a(x)=k_{a} \prod_{i=1}^N x_i^{a_i} \quad {\rm for}\  \ x\in B_1(0).
\end{equation}
Now we use the test function
$\zeta_{\varepsilon,a}(x):=\zeta_a(\varepsilon^{-1}x)$ for $ x\in\tilde\Omega$
in  (\ref{3.3}),
we have that
$$\sum_{|a|\le q} k_{a} D^{a}\zeta_{\varepsilon,a}(0)=\frac{k_{a}^2}{\varepsilon^{|a|}} \prod^{N}_{i=1}a_i! ,$$
where  $ a_i!=a_i\cdot (a_i-1)\cdots1>0$  and $ a_i!=1$ if $a_i=0$.

Let $r>0$, we obtain that
\begin{eqnarray*}
\Big|\int_{\tilde \Omega} \tilde u\mathcal{L}_\beta^* \zeta_{\varepsilon}\,d\tilde \gamma_\beta\Big|  &=& \Big|\int_{B_{2\varepsilon}(0)} \tilde u\mathcal{L}_\beta^* \zeta_{\varepsilon}\,d\tilde \gamma_\beta\Big|
\\&\le&\frac1{\varepsilon^{2}} \Big|\int_{B_{2\varepsilon}(0)} \tilde u(x)(-\Delta)   \zeta_a(\varepsilon^{-1} x)\,d\tilde \gamma_\beta\Big|
\\&&+\frac{2|\tau_+(\beta)|}{\varepsilon} \,\Big|\int_{B_{2\varepsilon}(0)}\tilde u(x)  \frac{x}{|x|^2}\cdot \nabla\zeta_a(\varepsilon^{-1} x)\,d\tilde \gamma_\beta\Big|
\\&\le&c_{13}  \left[ \frac1{\varepsilon^{2}}\int_{B_{2\varepsilon}(0)}| \tilde u(x)|\,d\tilde \gamma_\beta+\frac{1}{\varepsilon} \,\int_{B_{2\varepsilon}(0)} \frac{|\tilde u(x)|}{|x|}\,d\tilde \gamma_\beta\right]
\\&\le &\frac{c_{14}}{\varepsilon} \,\int_{B_{2\varepsilon}(0)} \frac{|\tilde u(x)|}{|x|}\,d\tilde \gamma_\beta ,
\end{eqnarray*}
then, by the fact that $u\in L^1(\Omega,\frac1{|x|}d\gamma_\beta)$, it follows  that
\begin{equation}\label{3.0.0}
\lim_{\varepsilon\to0^+}\D\int_{B_{2\varepsilon}(0)} \frac{|\tilde u(x)|}{|x|}\,d\tilde \gamma_\beta=0\quad{\rm and}\quad \lim_{\varepsilon\to0^+}\varepsilon \Big|\int_{\tilde \Omega}\tilde u\mathcal{L}_\beta^* \zeta_{\varepsilon}\,d\tilde \gamma_\beta \Big|=0.
\end{equation}
For $|a|\ge 1$, we have that
$$k_{a}^2\le c_{15}\varepsilon^{|a|-1}\Big|\int_{\tilde \Omega}\tilde u\mathcal{L}_\beta^* \zeta_{\varepsilon}\,d\tilde \gamma_\beta \Big|\to 0\quad{\rm as}\quad \varepsilon\to0,$$
then we have $k_{a}=0$ by arbitrary of $\varepsilon>0$ in (\ref{3.3}) with $|a|\ge1$, thus,
\begin{equation}\label{3.3-1}
L(\zeta)= \int_{\tilde \Omega} \left[\tilde u \mathcal{L}_\beta^*\zeta -\tilde f\zeta\right]\,d\tilde \gamma_\beta=k_0\zeta(0),\quad\  \forall\, \xi\in C^\infty_c (\tilde \Omega).
\end{equation}
For any $\zeta\in C^{1.1}_c(\tilde\Omega)$,  by taking a sequence $\zeta_n\in C^\infty_c (\tilde \Omega)$  converging to $\zeta$, we obtain that
(\ref{3.3-1}) holds for any $\zeta\in C^{1.1}_c(\tilde\Omega)$.

\vskip2mm
Now we fix  $\xi\in C^{1.1}_0(\Omega)$ with compact support in $\Omega\cup \{(x',0)\in \R^{N-1}\times\R: |x'|<r_0\}$, then $\xi/x_N\in C^{1.1}(\bar\Omega)$ and we may do  $x_N$-even extension of $\xi/x_N$ in $\tilde \Omega$, denoting by $\tilde \xi$,
then $\tilde \xi\in C^{1.1}_c(\tilde\Omega)$,  by the   $x_N$-even extension,
we have that
$$\tilde \xi(0)=\frac{\partial \xi}{\partial x_N}(0).$$
So it follows from (\ref{3.3-1}) that
\begin{equation}\label{3.3-2}
 \int_{\Omega} \Big(\tilde u \mathcal{L}_\beta^*(\frac{\xi}{x_N}) -\tilde f\frac{\xi}{x_N}\Big)\,d\gamma_\beta=k_0\frac{\partial \xi}{\partial x_N}(0),\quad\   \forall\, \xi\in C^\infty_c (\tilde \Omega),
\end{equation}
so (\ref{1.2f}) holds. \smallskip

\vskip2mm
$(ii)$ By the linearity of $\mathcal{L}_\beta$, we may assume that $f\ge0$. Let $f_n=f\eta_n$, where
$\eta_n(r)=1-\eta_0(nr)$ for $r\ge0$, where $\eta_0$ satisfies (\ref{eta}) and let $v_n$ be solution of
(\ref{2.02}) where  $f$ is replaced by $f_n$. We see that  $f_n$ is bounded and for any $\xi\in C^{1.1}_0(\Omega)$,
\begin{equation}\label{4.1-1}
 \int_{\Omega} v_n\,   \mathcal{L}_\beta^*(\frac{\xi}{x_N} )\, d\gamma_\beta =\int_{\Omega} f_n  \frac{\xi}{x_N} \, d\gamma_\beta.
\end{equation}
Then taking $\xi=w_2$ in Lemma \ref{lm 2.5}, we have that $ v_n$ is uniformly bounded in $L^1(\Omega, \,d\gamma_\beta)$  and  in $L^1(\Omega, \,x_N^{-1} d\gamma_\beta)$, that is,
$$\norm{v_n}_{L^1(\Omega, \,x_N^{-1} d\gamma_\beta)}\le \norm{\frac{\xi}{x_N}}_{L^\infty(\Omega)}\norm{f_n}_{L^1(\Omega, d\gamma_\beta)}\le \norm{\frac{\xi}{x_N}}_{L^\infty(\Omega)}\norm{f}_{L^1(\Omega, d\gamma_\beta)}. $$
Moreover,  $ \{v_n\}$ is increasing, and then there exists $v_f$ such that
$$v_n\to v_f\quad{\rm a.e.\ in}\ \ \Omega\quad{\rm and\ in}\ \ L^1(\Omega,\,x_N^{-1} d\gamma_\beta).$$
Then we have that
$$\int_{\Omega} v_f  \mathcal{L}_\beta^*(\xi )\, d\gamma_\beta =\int_{\Omega} f   \xi \,  d\gamma_\beta,\quad\forall\, \xi\in C^{1.1}_0(\Omega).$$
Since $f\in C^\gamma(\overline{\Omega}\setminus \{0\})$, then it follows by the standard regularity theory that
$v_f\in C^2(\Omega)$.

\vskip2mm
We claim that $v_f$ is a classical solution of (\ref{eq 1.1f}). From Corollary 2.8 in \cite{V} with $L^*=\mathcal{L}_{\beta}^*$, which is strictly elliptic in $\Omega\setminus B_r(0)$, we have that for $q<\frac{N}{N-1}$,
\begin{eqnarray}
\norm{v_n\lambda_\beta}_{W^{1,q}(\Omega_{2r})}  &\le & c_{16}\norm{f\lambda_\beta}_{L^1 (\Omega\setminus B_r(0))}+ c_{16}\norm{v_n\lambda_\beta}_{L^1(\Omega\setminus B_r(0))} \nonumber\\
  &\le & c_{17} \norm{f }_{L^1 (\Omega,\,d\gamma_\beta)},\label{5.1}
\end{eqnarray}
where $\Omega_{2r}=\{x\in\Omega\setminus B_{2r}(0):\, \rho(x)>2r\}.$
 We see that
\begin{eqnarray*}
-\Delta v_n   = -\frac{\beta}{|x|^2}v_n + f.
\end{eqnarray*}
For any compact set $K$ in $\Omega$, it is standard to improve the regularity $v_n$
$$\norm{v_n}_{C^{2,\lambda}(K)}\le c_{18}[\norm{f }_{L^1 (\Omega,\,d\gamma_\beta)}+ \norm{f}_{C^\lambda (K)}]$$
where $c_{18}>0$ is independent of $n$. Then  $v_f$ is a classical solution of  (\ref{eq 1.1f}) verifying the identity
 \begin{equation}\label{1.2f0}
 \int_{\Omega}v_f\,  \mathcal{L}_\beta^*(\frac{\xi}{x_N})\, d\gamma_\beta  = \int_{\Omega}   \frac{f\xi}{x_N} \, d\gamma_\beta,\quad\forall\, \xi\in   C^{1.1}_0(\Omega).
\end{equation}
 Together with  the fact that $u_{k,f}=k\Lambda^\Omega_\beta+v_f$, we conclude that the function $u_{k,f}$ is  a solution of (\ref{eq 1.1f}), verifying the identity (\ref{1.2f}) by (\ref{1.2f0}).

\vskip2mm
Finally, we prove the uniqueness. In fact, let $w_{k,f}$ be  a solution of (\ref{eq 1.1f})  verifying the identity (\ref{1.2f}).
$$\int_{\Omega}( u_{k,f}-w_{k,f})   \mathcal{L}_\beta^*(\frac{\xi}{x_N})\, d\gamma_\beta =0. $$
 For any Borel subset $O$  of $\Omega$, Corollary \ref{re 3.1} implies that problem
\begin{equation}\label{L00}
\arraycolsep=1pt\left\{
\begin{array}{lll}
 \displaystyle    \mathcal{L}_\beta^*(\frac{u}{x_N})   = \zeta_n\quad
   &{\rm in}\;\;   \Omega,\\[1mm]
 \phantom{  \mathcal{L}_\mu^*  }
 \displaystyle  u= 0\quad & {\rm   on}\;\; \partial{\Omega},
 \end{array}\right.
 \end{equation}
 has  a solution  $\eta_{\omega,n}\in C^2(\Omega)\cap C^{0.1}_0(\Omega)$,
where $\zeta_n:\bar\Omega\mapsto[0,1]$ is a $C^1(\bar\Omega)$
function such that $\zeta_n\to\chi_O\; {\rm{in}}\ L^\infty( \Omega)\; {\rm{as}}\ n\to\infty.$
Therefore by
passing to the limit as $n\to\infty$, we have that
$$\displaystyle\int_{O}( u_{k,f}-w_{k,f})d\gamma_\beta =0,$$
which implies that $u_{k,f}=w_{k,f}$ a.e. in $\Omega$ and then the uniqueness holds true. \hfill$\Box$

\begin{remark}\label{re 3.2}
Let $u_f$ be the solution of (\ref{eq 1.1f}) verifying the identity (\ref{1.2f}) with $k=0$, then
$u_f$ satisfies the isolated singular behavior (\ref{1.3}). In fact, letting $f\ge0$, then $u_f\ge 0$ in $\Omega$. So if (\ref{1.3}) fails, it implies by the positivity of
$u_f$, that $ \liminf_{t\to0^+}\inf_{z\in S^{N-1}_+}\frac{u_{f}(tz)}{\Lambda_\beta(tz)}=l_0>0$ and
$\tilde u_f:=u_f-l_0\Lambda_\beta^\Omega$ is a solution of (\ref{eq 1.1f}). By Lemma \ref{lm cp}, we have that $\tilde u_f\ge0$ in $\Omega$,
By the approximating procedure,
$\tilde u_f$ verifies the identity (\ref{1.2f}) with $k=0$, which is impossible with the fact that
$u_f-\tilde u_f=l_0\Lambda_\beta^\Omega$, which satisfies
$$ \int_{\Omega} (u_f-\tilde u_f) \mathcal{L}_\beta^*(\frac{\xi}{x_N} )\, d\gamma_\beta =l_0 c_\beta\frac{\partial \xi}{\partial x_N}(0),\quad \forall\,\xi\in  C^{1.1}_0(\Omega).$$

\end{remark}


\subsection{ Nonzero Dirichlet boundary}

Recall that  $P_\Omega$  is Poisson's Kernel of $-\Delta$ in $\Omega$   and
$ \mathbb{P}_\Omega[g](x)=\displaystyle\int_{\partial \Omega} P_\Omega(x,y)g(y)d\omega(y).$
It is known that if $g$ is continuous,  $\mathbb{P}_\Omega[g]$ is a solution of
\begin{equation}\label{6.1}
\arraycolsep=1pt\left\{
\begin{array}{lll}
 \displaystyle  -\Delta u = 0\quad
   {\rm in}\;\;  \Omega,\\[1mm]
 \phantom{ -\Delta }
 \displaystyle  u= g\quad  {\rm   on}\;\; \partial{\Omega}.
 \end{array}\right.
\end{equation}
Multiply $\frac{\xi\lambda_\beta}{x_N}$ where $\xi\in   C^{1.1}_0(\Omega)$ and integrate over $\Omega$, then we have that
\begin{eqnarray*}
0 &=& \int_\Omega (-\Delta\mathbb{P}_\Omega[g])\frac{\xi\lambda_\beta}{x_N} dx   \\
   &=& \int_{\partial\Omega}  \mathbb{P}_\Omega[g] \nabla(\frac{\xi\lambda_\beta}{x_N})\cdot\nu d\omega+ \int_\Omega \mathbb{P}_\Omega[g]
   \Big(-\Delta (\frac{\xi\lambda_\beta}{x_N}) \Big)dx \\
   &=&\int_{\partial\Omega} g \frac{\partial \xi}{\partial \nu} d\omega_\beta+\int_{\Omega}  \mathbb{P}_\Omega[g]  \mathcal{L}_\beta^*(\frac{\xi}{x_N})\, d\gamma_\beta-\beta\int_{\Omega}\frac{\mathbb{P}_\Omega[g] }{|x|^2}\frac\xi{x_N} d\gamma_\beta,
\end{eqnarray*}
that is, for any $\xi\in C^{1.1}_0(\Omega)$, there holds
 \begin{equation}\label{g}
 \int_{\Omega}  \mathbb{P}_\Omega[g]  \mathcal{L}_\beta^*(\frac{\xi}{x_N})\, d\gamma_\beta  =\beta\int_{\Omega}\frac{\mathbb{P}_\Omega[g] }{|x|^2}\frac\xi{x_N} d\gamma_\beta-\int_{\partial\Omega} g\frac{\partial \xi}{\partial \nu} d\omega_\beta.
\end{equation}

\begin{lemma}\label{lm 5.1}
Let $\beta\in[ \beta_0,\,+\infty)\setminus\{0\}$, $d\tilde\omega_\beta=(1+|x|^{\tau_+(\beta)})d\omega(x)$ and $g\ge0$. We have that

$(i)$ If  $g\in C(\partial\Omega\setminus\{0\})\cap L^1(\partial\Omega,\,d\tilde\omega_\beta)$, then
$\frac1{|\cdot \; |^2}\mathbb{P}_\Omega[g] \in   L^1(\Omega,\, d\gamma_\beta).$

$(ii)$ If  $g\in C(\partial\Omega\setminus\{0\})$ and
\begin{equation}\label{g1}
\lim_{r\to0^+}\int_{\partial\Omega\setminus B_r(0)} g \,d\tilde\omega_\beta=+\infty,
\end{equation}
 then
$$\lim_{r\to0^+} \int_{ \Omega\setminus B_r(0)} \frac1{|x|^2}\mathbb{P}_\Omega[g](x)d\gamma_\beta=+\infty.$$

\end{lemma}
\noindent {\bf Proof.} From Proposition 2.1 in \cite{BV} that
\begin{equation}\label{3.00}
 c_{19}\rho(x)|x-y|^{-N}\le P_\Omega(x,y)\le c_{20}\rho(x)|x-y|^{-N},\quad x\in\Omega,\ y\in\partial\Omega,
\end{equation}
where $\rho(x)={\rm dist}(x,\partial\Omega)$. Since $g$ is continuous in $\partial\Omega\setminus\{0\}$ and $\Omega$ is flat near the origin, we can only consider the integrability of $\frac1{|\cdot \;|^2}\mathbb{P}_\Omega[g]$ near the origin. Fix $r=r_0/2$,
let $B'_{r}(0)=\{x' \in\R^{N-1}: |x'|<r\}$ and $e_{(y',0)}=(\frac{y'}{|y'|},\, 0)$ for $y'\not=0$,  then
\begin{eqnarray*}
&& \int_{B_{r}^+(0)} \frac1{|x|^2}\mathbb{P}_\Omega[g] d\gamma_\beta
\\&\ge & c_{21} \int_{B_{r}^+(0)} \int_{B'_{r}(0)\setminus\{0\}} g(y') |x-(y',0)|^{-N}\frac{x_N^2}{|x|^2}  |x|^{\tau_+(\beta)}\, dy'  dx \\
   &=& c_{22}\int_{B_{r}'(0)\setminus\{0\}}g(y') |y'|^{\tau_+(\beta)} \int_{B_{r/|y'|}^+(0)}    |z-e_{(y',0)}|^{-N}\frac{z_N^2}{|z|^2}  |z|^{\tau_+(\beta)} dz   dy'  
\end{eqnarray*}
and
\begin{eqnarray*}
&& \int_{B_{r}^+(0)} \frac1{|x|^2}\mathbb{P}_\Omega[g] d\gamma_\beta
\\&\le & c_{23} \int_{B_{r}^+(0)} \int_{B'_{r}(0)\setminus\{0\}} g(y') |x-(y',0)|^{-N}\frac{x_N^2}{|x|^2}  |x|^{\tau_+(\beta)}\, dy'  dx \\
   &=& c_{24}\int_{B_{r}'(0)\setminus\{0\}}g(y') |y'|^{\tau_+(\beta)} \int_{B_{r/|y'|}^+(0)}    |z-e_{(y',0)}|^{-N}\frac{z_N^2}{|z|^2}  |z|^{\tau_+(\beta)} dz   dy'.  
\end{eqnarray*}
Now we do estimates for $$\int_{B_{r/|y'|}^+(0)}  I(z) dz:=\int_{B_{r/|y'|}^+(0)}    |z-e_{(y',0)}|^{-N}\frac{z_N^2}{|z|^2}  |z|^{\tau_+(\beta)} dz,$$
we have
\begin{eqnarray*}
0<\int_{B_{\frac12}^+(0)}  I(z)\, dz     \le  2^{N}\int_{B_{\frac12}^+(0)} |z|^{\tau_+(\beta)} dz,
\end{eqnarray*}
\begin{eqnarray*}
0<\int_{B_{\frac12}^+(e_{(y',0)})}  I(z)\, dz    & \le & 2^{|\tau_+(\beta)|+2}\int_{B_{\frac12}^+(e_{(y',0)})}   |z-e_{(y',0)}|^{-N} z_N^2 dz
\\&\le &2^{|\tau_+(\beta)|+2}\int_{B_{\frac12}^+(0)}   |z|^{2-N}  dz,
\end{eqnarray*}
and
\begin{eqnarray*}
&&\int_{B_{r/|y'|}^+(0)\setminus\left(B_{\frac12}^+(0)\cup B_{\frac12}^+(e_{(y',0)})\right)}  I(z) dz
\\  & \le & c_{25}\int_{B_{r/|y'|}^+(0)\setminus B_{\frac12}^+(0)}  |z|^{-N+\tau_+(\beta)}\,dz
\\&\le & \left\{
\begin{array}{lll}
 \displaystyle  c_{26}\int_{\R^N \setminus B_{\frac12}(0)}  |z|^{-N+\tau_+(\beta)}\,dz  \qquad
   &{\rm if}\quad  \beta<0,\\[2mm]
 \phantom{   }
 \displaystyle  c_{26}|y'|^{-\tau_+(\beta)} \qquad  &{\rm   if}\quad \beta>0
 \end{array}\right.
 \\&\le &c_{27}(1+|y'|^{-\tau_+(\beta)})
\end{eqnarray*}
and
\begin{eqnarray*}
\int_{B_{r/|y'|}^+(0)\setminus\left(B_{\frac12}^+(0)\cup B_{\frac12}^+(e_{(y',0)})\right)}  I(z)\, dz    & \ge & c_{28}\int_{B_{r/|y'|}^+(0)\setminus B_{\frac12}^+(0)}  |z|^{-N+\tau_+(\beta)}\,dz
\\&\ge & c_{29}(1+|y'|^{-\tau_+(\beta)}).
\end{eqnarray*}
Thus, we have that
\begin{equation}\label{5.2}
c_{30}\int_{B_{r}'(0)\setminus\{0\}}g(y')  d\tilde\omega(y') \le \int_{B_{r}^+(0)} \frac1{|x|^2}\mathbb{P}_\Omega[g] d\gamma_\beta \le c_{31}\int_{B_{r}'(0)\setminus\{0\}}g(y') d\tilde\omega(y'),
\end{equation}
which,  together with  the fact that $\mathbb{P}_\Omega[g] $ is nonnegative  and bounded in $\Omega\setminus B_{r}^+(0)$,
proves  Lemma \ref{lm 5.1}.\hfill$\Box$\medskip

We remark that Lemma \ref{lm 5.1} provides estimates for transforming  the boundary data into the nonhomogeneous term. Now we are ready to prove
Theorem \ref{teo 1} part $(i)$ where we distinguish two cases $\beta\in[\beta_0,\, 0]$ and $\beta>  0$.

\vskip3mm
\noindent{\bf Proof of Theorem \ref{teo 1}.  Part $(i)$.}
{\it The existence for $g\in L^1(\partial\Omega,\, d\tilde\omega_\beta)$. }
Let $\bar f=f-\frac{\beta}{|\cdot \;|^2}\mathbb{P}_\Omega[g].$  Then it follows from Lemma \ref{lm 5.1} part $(i)$ that $\bar f\in L^1(\Omega,\,d\gamma_\beta)$ and applying Theorem \ref{teo 3.1} part $(i)$,  problem (\ref{eq 1.1f}) verifying (\ref{1.2f}) for  $k\in\R$ and replaced $f$  by $\bar f$  admits a unique solution of  $u_f$.
 Denote $u_{f,g}:=u_f+\mathbb{P}_\Omega[g]$,
then
$$\mathcal{L}_\beta u_{f,g}=f\quad{\rm and}\quad u_{f,g}=g\quad {\rm  on}\ \ \partial\Omega\setminus\{0\}.$$
Together with (\ref{1.2f}) and (\ref{g}), we have that
$u_{f,g}$ verifies (\ref{1.2fg})
and it is the unique solution of  problem (\ref{eq 1.1f}) verifying (\ref{1.2f}) for that $k$. \smallskip

 {\it Case of $\beta\in[\beta_0,\, 0]$}. Then $d\tilde\omega_\beta$ is equivalent to $d\omega_\beta$, so  $L^1(\partial\Omega,\, d\tilde \omega_\beta)= L^1(\partial\Omega,\, d\omega_\beta)$
and  we are done.
 \smallskip

  {\it Case of $\beta>0$. }    We note that
$$L^1(\partial\Omega,\, d\tilde \omega_\beta)\subsetneqq  L^1(\partial\Omega,\, d\omega_\beta).$$
So for $g\in L^1(\partial\Omega,\, d\omega_\beta)\setminus L^1(\partial\Omega,\, d\tilde\omega_\beta)$, we may assume $g\ge0$ by linearity of $\mathcal{L}_\beta$. Let
\begin{equation}\label{6-1}
 \eta_n(s)=1-\eta_0(ns)\quad{\rm and}\quad g_n(x)=g(x)\eta_n(|x|),
\end{equation}
where $\eta_0$ is defined in (\ref{eta}). Then $\{g_n\}_n \subset L^1(\partial\Omega,\,d\tilde \omega_\beta)$ is an increasing  sequence of functions.
For simplicity, we assume that  $f=0$.  Then the problem
\begin{equation}\label{eq 5.1.1}
\left\{\arraycolsep=1pt
\begin{array}{lll}
\mathcal{L}_\beta^* u= 0\quad &{\rm in}\;\;  \Omega, \\[1mm]
 \phantom{ \mathcal{L}_\beta  }
  u=g_n\quad &{\rm on}\;\; \partial\Omega\setminus\{0\}
  \end{array}
 \right.
\end{equation}
has a unique solution of $u_{n}$ verifying
the identify
 \begin{equation}\label{1.2fg0}
 \int_{\Omega}u_{n} \mathcal{L}_\beta^*(\frac{\xi}{x_N})\, d\gamma_\beta  =  -\int_{\partial\Omega} g_n\frac{\partial \xi}{\partial \nu} d\omega_\beta, \quad\forall\, \xi\in   C^{1.1}_0(\Omega).
\end{equation}
Since $0\le g_n\le g$ and $g\in  L^1(\partial\Omega,\, d\omega_\beta)$, we may expand the text function space  including   $w_1$, $w_2$, which are the solutions of (\ref{eq 2.2}) and (\ref{eq 2.2-1}) respectively.
 Taking $\xi=w_1$ and then $w_2$ ,   we derive that
$$\norm{u_n}_{L^1(\Omega)}\le c_{32}\norm{g_n}_{L^1(\partial\Omega,\, d\omega_\beta)}\le c_{33}\norm{g}_{L^1(\partial\Omega,\, d\omega_\beta)}$$
and
$$\norm{u_n}_{L^1(\Omega, \,x_N^{-1} d\gamma_\beta)}\le c_{34}\norm{g}_{L^1(\partial\Omega,\, d\omega_\beta)}.$$

We notice that $u_{n}\ge0$ and  the mapping $n\mapsto u_{n}$ is increasing, then by the monotone converge theorem,
we have that there exists $u$ such that $u_n$ converging to $u$ in $L^1(\Omega, \frac1{x_N} d\gamma_\beta)$.
Since $\xi\in C_0^{1.1}(\Omega)$, we have that $|\mathcal{L}_\beta^*(\xi/x_N)|\le cx_N^{-1}.$
Pass to the limit of (\ref{1.2fg0}), we have that $u$ verifies  that
 \begin{equation}\label{1.2fg00}
 \int_{\Omega}u \mathcal{L}_\beta^*(\xi/x_N)\, d\gamma_\beta  =  -\int_{\partial\Omega} g\frac{\partial \xi}{\partial \nu} d\omega_\beta, \quad\forall\, \xi\in   C^{1.1}_0(\Omega).
\end{equation}
From standard interior regularity, we have that $u$ is a classical solution
$$\left\{\arraycolsep=1pt
\begin{array}{lll}
\mathcal{L}_\beta^* u= 0\quad {\rm in}\;\;  \Omega, \\[1mm]
 \phantom{ \mathcal{L}_\beta  }
  u=g\quad {\rm on}\;\;  \partial\Omega\setminus\{0\},
  \end{array}
 \right.$$
which ends the proof.\hfill$\Box$

\setcounter{equation}{0}
\section{Nonexistence}

  In this subsection, we establish the approximation of the fundamental solution $G_\mu$.

\begin{lemma}\label{lm 2.3}

$(i)$ Let $\{\delta_n\}_n$ be a sequence  of nonnegative  $L^\infty$-functions defined in $\Omega$ such that
${\rm supp}\,\delta_n\subset B_{r_n}(0)\cap \Omega,$
 where $r_n\to0$ as $n\to+\infty$ and
$$ \int_\Omega \delta_n \xi dx \to \frac{\partial\xi(0)}{\partial x_N}\quad{\rm as} \quad n\to+\infty,\quad \forall \xi\in C_0^1(\Omega).$$
 For any $n$, let $w_n$  be the unique solution of the problem in the $d\gamma_\beta$-distributional sense
 \begin{equation}\label{eq 2.3}
 \arraycolsep=1pt\left\{
\begin{array}{lll}
 \displaystyle  \mathcal{L}_\beta  u = \delta_n/\lambda_\beta\qquad
  & {\rm in}\;\;  {\Omega}\setminus \{0\},\\[2mm]
 \phantom{  L_\beta }
 \displaystyle  u= 0\qquad  &{\rm   on }\;\;  \partial{\Omega},\\[2mm]
 \phantom{   }
  \displaystyle \lim_{r\to0^+} \sup_{x\in \partial_+ B_r(0)} \frac{|u(x)|}{\Lambda_\beta (x)} =0.
 \end{array}\right.
\end{equation}
Then
$$\lim_{n\to+\infty} w_n(x)=\frac1{c_\beta}\Lambda_\beta^\Omega(x),\quad \forall \, x\in\Omega\setminus\{0\}$$
and for any compact set $K\subset \Omega\setminus\{0\}$,
 \begin{equation}\label{3.7}
  w_n\to \frac1{c_\beta}\Lambda_\beta^\Omega  \quad{\rm as}\quad n\to+\infty \ \ {\rm in}\ \ C^{2}(K).
 \end{equation}

$(ii)$ Let $\{\sigma_n\}_n$ be a sequence  of nonnegative  $L^\infty$ functions defined on $\partial \Omega$ such that
${\rm supp}\,\sigma_n\subset \partial \Omega\cap B_{r_n}(0),$
 where $r_n\to0$ as $n\to+\infty$ and
$$ \int_{\partial \Omega} \sigma_n \zeta d\omega(x) \to   \zeta(0) \quad{\rm as} \quad n\to+\infty,\quad \forall \zeta\in C^1(\partial\Omega).$$
 For any $n$, let $v_n$  be the unique solution of the problem
 \begin{equation}\label{eq 2.3g}
 \arraycolsep=1pt\left\{
\begin{array}{lll}
 \displaystyle  \mathcal{L}_\beta  u =0\quad
  & {\rm in}\;\; {\Omega}\setminus \{0\},\\[2mm]
 \phantom{   }
 \displaystyle  u= \frac{\sigma_n}{|\cdot|^{\tau_+(\beta)}} \quad  &{\rm   on }\;\;  \partial{\Omega}\setminus\{0\}
 \end{array}\right.
\end{equation}
subject to
$$\int_{\Omega}v_n \mathcal{L}_\beta^*(\xi/x_N)\, d\gamma_\beta  = -\int_{\partial\Omega} \sigma_n\frac{\partial \xi}{\partial \nu} d\omega,\quad\forall\, \xi\in   C^{1.1}_0(\Omega).$$
Then
$$\lim_{n\to+\infty} v_n(x)=\frac1{c_\beta}\Lambda_\beta^\Omega(x),\quad \forall \, x\in\Omega\setminus\{0\}$$
and for any compact set $K\subset \Omega\setminus\{0\}$, (\ref{3.7}) holds true.

\end{lemma}
\noindent{\bf Proof.} From Lemma \ref{lm 2.1}, problems (\ref{eq 2.3}) and (\ref{eq 2.3g}) have unique solutions $w_n,v_n\ge0$   respectively and satisfying that
\begin{equation}\label{3.0-1}
  \int_{\Omega} w_n   \mathcal{L}_\beta^*(\frac{\xi}{x_N})\, d\gamma_\beta  =\int_{\Omega} \delta_n\xi\,  dx,\quad \forall\,\xi\in C^{1,1}_0(\Omega)
\end{equation}
and
\begin{equation}\label{3.0-2}
  \int_{\Omega} v_n  \mathcal{L}_\beta^*(\frac{\xi}{x_N})\, d\gamma_\beta =-\int_{\partial\Omega} \sigma_n\frac{\partial \xi}{\partial \nu} d\omega_\beta,\quad \forall\,\xi\in C^{1,1}_0(\Omega).
\end{equation}
By taking $\xi=\xi_0$, the solution of (\ref{eq 2.2}), we obtain that
$$\norm{w_n}_{L^1(\Omega,\,  d\gamma_\beta)}\le \norm{\xi_0}_{L^\infty(\Omega)}\norm{\delta_n}_{L^1(\Omega)}=\norm{\xi_0}_{L^\infty(\Omega)}. $$
For any  $r>0$, take $\xi$ with the support in $\Omega\setminus B_r(0)$,
then  $\xi\in C^{1.1}_c(\overline{\Omega\setminus B_r(0)})$,
$$\int_{\Omega\setminus B_r(0)} w_n  \mathcal{L}_{\mu}^* (\xi) \, d\gamma_\beta =0. $$
Take $\xi$ the solution of (\ref{eq 2.2f}) with $f(x)=\frac1{|x|}$, we have that
\begin{equation}\label{3.0-1.1}
  \int_{\Omega} w_n  |x|^{-1}\, d\gamma_\beta  =\int_{\Omega} \delta_n\xi\,  dx\le \norm{\xi_0}_{L^\infty(\Omega)}
\end{equation}
and
\begin{equation}\label{3.0-2.2}
  \int_{\Omega} v_n  |x|^{-1} \, d\gamma_\beta =-\int_{\partial\Omega} \sigma_n\frac{\partial \xi}{\partial \nu} d\omega_\beta\le \norm{\nabla \xi_0}_{L^\infty(\Omega)}.
\end{equation}
So  $w_n,v_n$ are uniform bounded in $L^1(\Omega, \,|x|^{-1} d\gamma_\beta)$.

From Corollary 2.8 in \cite{V} with $L^*=\mathcal{L}_{\mu}^*$, which is strictly elliptic in $\Omega\setminus B_r(0)$, we have that for $q<\frac{N}{N-1}$,
\begin{eqnarray*}
\norm{w_n\lambda_\beta}_{W^{1,q}(\Omega_{2r})}  \le  c_{35}\norm{\delta_n}_{L^1 (\Omega\setminus B_r(0))}+ c_{36}\norm{w_n }_{L^1(\Omega\setminus B_r(0),\,d\gamma_\beta)} \le c_{37}
\end{eqnarray*}
and
\begin{eqnarray*}
\norm{v_n\lambda_\beta}_{W^{1,q}(\Omega_{2r})}  &\le & c_{38}\norm{\sigma_n}_{L^1 (\partial\Omega\setminus B_r(0))}+ c_{39}\norm{v_n }_{L^1(\Omega\setminus B_r(0),\,d\gamma_\beta)}\le c_{40},
\end{eqnarray*}
where $\Omega_{2r}=\{x\in\Omega\setminus B_{2r}(0):\, \rho(x)>2r\}.$
By the compact embedding $W^{1,q}(\Omega_{2r})\hookrightarrow L^1(\Omega_{2r}),$
up to some subsequence, there exists $w_\infty,\,v_\infty\in W^{1,q}_{loc}(\Omega)\cap L^1(\Omega,\, d\gamma_\beta)$ such that
$$w_n \to w_\infty  \quad{\rm as}\quad n\to+\infty\quad {\rm a.e.\ \ in}\ \ \Omega\ \ {\rm and\ in}\quad L^1(\Omega,\,d\gamma_\beta)$$
and it follows by (\ref{3.0-1}) and (\ref{3.0-2}) that for $\xi\in C^{1.1}_0(\Omega)$,
$$\int_{\Omega} w_\infty   \mathcal{L}_\beta^*(\xi)\, d\gamma_\beta =\int_{\Omega} v_\infty   \mathcal{L}_\beta^*(\xi)\, d\gamma_\beta= \frac{\partial\xi}{\partial x_N}(0). $$
Furthermore,
$$\int_{\Omega}( w_\infty-\frac1{c_\beta}G_\beta)   \mathcal{L}_\beta^*(\xi)\, d\gamma_\beta =0. $$
From the Kato's inequality,  we deduce that
$$w_\infty=v_\infty=\frac1{c_\beta}\Lambda_\beta^\Omega\quad{a.e.}\;\; \Omega. $$

\noindent{\it Proof of (\ref{3.7}).} For any $x_0\in \Omega\setminus\{0\}$, let $r_0=\frac14\{|x_0|,\, \rho(x_0)\}$ and
$\mu_n=w_n\eta,$
where $\eta(x)=\eta_0(\frac{|x-x_0|}{r_0})$. There exists $n_0>0$ such that for $n\ge n_0$,
${\rm supp}\mu_n\cap B_{r_n}(0)=\emptyset.$
Then
\begin{eqnarray*}
-\Delta \mu_n (x) &=& -\Delta w_n (x)\eta(x)-2 \nabla w_n\cdot\nabla\eta-w_n \Delta\eta \\
   &=&  -2 \nabla w_n\cdot\nabla\eta-w_n \Delta\eta,
\end{eqnarray*}
where $\nabla\eta$ and $\Delta\eta$ are smooth.

We observe that $w_n\in W^{1,q}(B_{2r_0}(x_0))$ and
$-2\nabla w_n\cdot\nabla\eta-w_n \Delta\eta\in L^q(B_{2r_0}(x_0)),$
then we have that
$$\norm{\mu_n}_{ W^{2,q}(B_{r_0}(x_0))}\le c\norm{w_n}_{L^1(\Omega,\,  d\gamma_\beta)},$$
where $c>0$ is independent of $n$.
Thus,  $-2\nabla w_n\cdot\nabla\eta-w_n \Delta\eta\in W^{1,q}(B_{r_0}(x_0)),$
repeat above process $N_0$ steps, for $N_0$ large enough, we deduce that
$$\norm{w_n}_{C^{2,\gamma}(B_{\frac{r_0}{2^{N_0}}}(x_0))}\le c\norm{w_n}_{L^1(\Omega,\,  d\gamma_\beta)},$$
where $\gamma\in(0,1)$ and $c>0$ is independent of $n$.
As a conclusion, (\ref{3.7}) follows by Arzel\`{a}-Ascola theorem and Heine-Borel theorem.
The above process also holds for $v_n$. This ends the proof. \hfill$\Box$\medskip

\vskip2mm
\noindent{\bf Proof of Theorem \ref{teo 1}. Part $(ii)$.}
From (\ref{1.3fg-}), one of the following two cases holds true,
$${\rm case}\ 1: \lim_{r\to0^+}  \int_{\Omega\setminus B_r(0)} f \,d\gamma_\beta  = +\infty, \;\; {\rm  or\  case}\ 2: \;\; \lim_{r\to0^+} \int_{\partial\Omega\setminus B_r(0)} g \,d\omega_\beta =+\infty.$$

{\it Case 1.}   We argue by contradiction. Assume that problem (\ref{eq 1.1fg}) has a nonnegative solution of $u_f$.
Let $\{r_n\}_n$ be a sequence of strictly decreasing positive numbers converging to $0$.
From  the fact $f\in C_{loc}^\gamma(\overline{\Omega}\setminus \{0\})$, for any $r_n$ fixed, we have that
$$
 \lim_{r\to0^+} \int_{(B_{r_n}(0)\setminus B_r(0))\cap \Omega}f(x) d\gamma_\beta  =+\infty,
$$
then there exists $R_n\in (0,r_n)$ such that
$$
  \int_{(B_{r_n}(0)\setminus B_{R_n}(0))\cap \Omega}f d\gamma_\beta =n.
$$
Let $\delta_n=\frac1n \lambda_\beta f\chi_{B_{r_n}(0)\setminus B_{R_n}(0)}$, then the problem
$$
 \arraycolsep=1pt\left\{
\begin{array}{lll}
 \displaystyle  \mathcal{L}_\mu u\cdot\lambda_\beta= \delta_n\qquad
   {\rm in}\quad  {\Omega}\setminus \{0\},\\[1mm]
 \phantom{  L_\mu -- }
 \displaystyle  u= 0\qquad  {\rm   on}\quad \partial{\Omega},\\[1mm]
 \phantom{   }
  \displaystyle \lim_{x\to0}u(x)\Phi_\mu^{-1}(x)=0
 \end{array}\right.
$$
has a unique positive solution  $w_n$  satisfying (in the usual sense)
$$\int_{\Omega} w_n \mathcal{L}_\mu(\lambda_\beta\xi) dx=\int_{\Omega} \delta_n \xi dx,\quad\forall\, \xi\in C^{1.1}_0(\Omega).$$
For any $\xi\in C^{1.1}_0(\Omega)$, we have that
$$\int_\Omega w_n  \mathcal{L}_\mu^*(\xi)\,  d\gamma_\beta =\int_{\Omega} \delta_n \xi\, dx\to \frac{\partial\xi}{\partial x_N}(0)\quad{\rm as}\quad n\to+\infty. $$
Therefore, by Lemma \ref{lm 2.3}  for any compact set $\mathcal{K}\subset \Omega\setminus \{0\}$
$$\norm{w_n-\Lambda_\beta^\Omega}_{C^1(\mathcal{K})}\to 0\quad{\rm as}\quad  {n\to+\infty}.$$
We fix a point $x_0\in \Omega$ and let $r_0=\frac{1}{2} \min\{|x_0|,\, \rho(x_0)\}$ and $\mathcal{K}=  \overline{B_{r_0}(x_0)}$, then
there exists $n_0>0$ such that for $n\ge n_0$,
\begin{equation}\label{4.1}
 w_n\ge \frac12G_\mu\quad{\rm in}\quad \mathcal{K}.
\end{equation}

Let $u_n$ be the solution (in the usual sense) of
$$
 \arraycolsep=1pt\left\{
\begin{array}{lll}
 \displaystyle  \mathcal{L}_\mu u\cdot\lambda_\beta= n\delta_n\quad
   {\rm in}\;\;  {\Omega}\setminus \{0\},\\[1mm]
 \phantom{  L_\mu -- }
 \displaystyle  u= 0\quad \ \ {\rm   on}\;\; \partial{\Omega},\\[1mm]
 \phantom{   }
  \displaystyle\lim_{r\to0^+} \sup_{x\in \partial_+ B_r(0)} \frac{|u(x)|}{\Lambda_\beta (x)} =0,
 \end{array}\right.
$$
then we have that $u_n\ge nw_n\;\; {\rm in}\;\;  \Omega.$
Together with (\ref{4.1}), we derive that
$$u_n\ge  \frac n2\Lambda_\mu^\Omega\quad {\rm in}\;\;  \mathcal{K}.$$
Then by comparison principle, we have that $u_f(x_0)\ge u_n(x_0)\to+\infty\quad{\rm as}\;\; n\to+\infty,$
which contradicts to the fact that $u_f$ is classical solution of (\ref{eq 1.1f}). \medskip

 {\it Case 2.} Similarly for any $n\in\N$, we can take
  $r_n>R_n>0$ such that $r_n\to0$ as $n\to+\infty$ and
  $$
  \int_{(B_{r_n}(0)\setminus B_{R_n}(0))\cap \partial \Omega}g d\omega_\beta =n.
$$
 Let $\sigma_n=\frac1n   g\chi_{B_{r_n}(0)\setminus B_{R_n}(0)} $,  $w_n$ be the solution of
 $$
 \arraycolsep=1pt\left\{
\begin{array}{lll}
 \displaystyle  \mathcal{L}_\mu u =0 \quad
   {\rm in}\;\;  {\Omega}\setminus \{0\},\\[1mm]
 \phantom{  L_\mu   }
 \displaystyle  u= \sigma_n/|\cdot|^{\tau^+(\beta)}\quad \ \ {\rm   on}\;\; \partial{\Omega},
 \end{array}\right.
$$
subject to
$$\int_{\Omega}w_n \mathcal{L}_\beta^*(\frac{\xi}{x_N})\, d\gamma_\beta  = -\int_{\partial\Omega} \sigma_n\frac{\partial \xi}{\partial \nu} d\omega,\quad\forall\, \xi\in   C^{1.1}_0(\Omega).$$
Repeat the procedure in Case 1, we get a contradiction which completes the proof.
  \hfill$\Box$

\bigskip

\bigskip

{\small \noindent{\bf Acknowledgements:}
H. Chen is supported by the Natural Science Foundation of China [Nos. 11661045, 11726614].
A. Quaas is partially supported by  Fondecyt Grant No. 1151180, Programa Basal, CMM in U. de Chile and Millennium Nucleus Center for Analysis of PDE NC130017.
F. Zhou is partially supported by NSFC  [Nos. 11726613, 11431005]; and STCSM [No. 18dz2271000].}

\end{document}